# ASYMPTOTICS FOR SLICED AVERAGE VARIANCE ESTIMATION[1]

By Yingxing Li and Li-Xing Zhu

*Cornell University and Hong Kong Baptist University*

In this paper, we systematically study the consistency of sliced average variance estimation (SAVE). The findings reveal that when the response is continuous, the asymptotic behavior of SAVE is rather different from that of sliced inverse regression (SIR). SIR can achieve $\sqrt{n}$ consistency even when each slice contains only two data points. However, SAVE cannot be $\sqrt{n}$ consistent and it even turns out to be *not consistent* when each slice contains a fixed number of data points that do not depend on $n$, where $n$ is the sample size. These results theoretically confirm the notion that SAVE is more sensitive to the number of slices than SIR. Taking this into account, a bias correction is recommended in order to allow SAVE to be $\sqrt{n}$ consistent. In contrast, when the response is discrete and takes finite values, $\sqrt{n}$ consistency can be achieved. Therefore, an approximation through discretization, which is commonly used in practice, is studied. A simulation study is carried out for the purposes of illustration.

## 1. Introduction.

Dimension reduction has become one of the most important issues in regression analysis because of its importance in dealing with problems with high-dimensional data. Let $Y$ and $\boldsymbol{x} = (x_1, \ldots, x_p)^T$ be the response and $p$-dimensional covariate, respectively. In the literature, when $Y$ depends on $\boldsymbol{x} = (x_1, \ldots, x_p)^T$ through a few linear combinations $B^T\boldsymbol{x}$ of $\boldsymbol{x}$, where $B = (\beta_1, \ldots, \beta_k)$, there are several proposed methods for estimating the projection directions $B$/space that is spanned by $B$, such as projection pursuit regression (PPR) [11], the alternating conditional expectation (ACE) method [1], principal Hessian directions (pHd) [17], minimum average variance estimation (MAVE) [23], iterated pHd [7] and profile least-squares

Received February 2005; revised February 2006.
[1]Supported by Grant HKU 7058/05P from the Research Grants Council of the Hong Kong SAR government, Hong Kong, China.

estimation [10]. All of these methods estimate the projection directions $B$ or the subspace that is spanned by $B$ when $B$ is contained within the mean regression function.

For more general models in which some $\beta_i$ are in the variance component of the model, two estimation methods—sliced inverse regression (SIR) [16] and sliced average variance estimation (SAVE) [5, 9]—have received much attention. SIR is based on the estimation of the conditional mean and SAVE on the estimation of the conditional variance function of the covariates given the response, the inverse regression. The aim of these two methods is to estimate the *central dimension reduction* (CDR) *space* that is defined as follows. Suppose that $Y$ is independent of $\boldsymbol{x}$, given $B^T\boldsymbol{x}$, which is written as $Y \perp\!\!\!\perp \boldsymbol{x}|B^T\boldsymbol{x}$, where $\perp\!\!\!\perp$ stands for independence and $B = (\beta_1, \dots, \beta_k)$ is an unknown $p \times k$ matrix, the columns of which are of unit length under the Euclidean norm and mutually orthogonal. A *dimension reduction subspace* is defined as the space that is spanned by the column vectors of $B$ and a CDR subspace is the intersection of all of the dimension reduction subspaces that satisfy conditional independence (see [3, 4]). The CDR subspace is still a dimension reduction subspace with the notation $S_{y|\boldsymbol{x}}$ under certain regularity conditions. SIR and SAVE are used to estimate $S_{y|\boldsymbol{x}}$. If we let $\boldsymbol{z} = \Sigma_{\boldsymbol{x}}^{-1/2}(\boldsymbol{x} - \mathbf{E}(\boldsymbol{x}))$ be the standardized covariate, then $S_{y|\boldsymbol{z}} = \Sigma_{\boldsymbol{x}}^{1/2} S_{y|\boldsymbol{x}}$ (see [4] for details). Hence, the estimation can be carried out equivalently for the pair of variables $(y, \boldsymbol{z})$. For convenience, we first use the standardized variable $\boldsymbol{z}$ to study the asymptotic behavior. In practice, the sample covariance matrix and the sample mean must be estimated and thus the results involving the estimated covariate $\hat{\boldsymbol{z}} = \hat{\Sigma}_{\boldsymbol{x}}^{-1/2}(\boldsymbol{x} - \bar{\boldsymbol{x}})$ will be reported as corollaries, where $\hat{\Sigma}_{\boldsymbol{x}}$ and $\bar{\boldsymbol{x}}$ are the sample covariance matrix and sample mean of the $\boldsymbol{x}_i$'s, respectively.

Denote the inverse regression function by $\mathbf{E}(\boldsymbol{z}|Y = y)$ and the conditional covariance of $\boldsymbol{z}$ given $y$ by $\Sigma_{\boldsymbol{z}|y} := \mathbf{E}((\boldsymbol{z} - \mathbf{E}(\boldsymbol{z}|Y))(\boldsymbol{z} - \mathbf{E}(\boldsymbol{z}|Y))^T|Y = y)$. SIR estimates the CDR subspace via the eigenvectors that are associated with the nonzero eigenvalues of the covariance matrix $\mathrm{Cov}(\mathbf{E}(\boldsymbol{z}|Y))$; SAVE estimates it via the eigenvectors that are associated with the nonzero eigenvalues of the covariance matrix $\mathbf{E}((I_p - \Sigma_{\boldsymbol{z}|Y})(I_p - \Sigma_{\boldsymbol{z}|Y})^T)$. For SIR estimation, we need the linearity condition

$$\mathbf{E}(\boldsymbol{z}|P_{S_{y|\boldsymbol{z}}}\boldsymbol{z}) = P_{S_{y|\boldsymbol{z}}}\boldsymbol{z}. \tag{1.1}$$

For SAVE estimation we also assume that

$$\mathrm{Cov}(\boldsymbol{z}|P_{S_{y|\boldsymbol{z}}}\boldsymbol{z}) = I_p - P_{S_{y|\boldsymbol{z}}}, \tag{1.2}$$

where $P_{(\cdot)}$ stands for the projection operator with respect to the standard inner product.

It is worth pointing out that the study of SAVE should receive more attention, as several papers have revealed that SAVE is more comprehensive than SIR: under regularity conditions, the CDR space of SAVE actually contains that of SIR (see [6, 24]). In particular, SIR will fail to work in symmetric regressions with $y = f(B^T\boldsymbol{x}) + \varepsilon$, where $f$ is a symmetric function of the argument $B^T\boldsymbol{x}$. Therefore, theoretically, SAVE should be a more powerful method than SIR under regularity conditions to estimate the CDR space.

Clearly, the primary aim is to estimate either $\mathrm{Cov}(\mathbf{E}(\boldsymbol{z}|Y))$ or $\mathbf{E}[(I_p - \Sigma_{\boldsymbol{z}|Y})(I_p - \Sigma_{\boldsymbol{z}|Y})^T]$. Li [16] proposed a slicing estimation that involves a very simple and easily implemented algorithm to estimate the inverse regression function, in which the slicing estimator is the weighted sum of the sample covariances of $\boldsymbol{z}_i$'s in each slice of $y_i$'s. He also demonstrated, by means of a simulation, that the performance of the slicing estimator is not sensitive to the choice of the number of slices. Zhu and Ng [27] provided a theoretical background for Li's empirical study and proved that $\sqrt{n}$ consistency and asymptotic normality hold provided the number of slices is within the range $\sqrt{n}$ to $n/2$. In other words, $\sqrt{n}$ consistency can be ensured when each slice contains a number of points between 2 and $\sqrt{n}$. The only thing that is affected by different numbers of slices is the asymptotic variance of the estimator. A relevant reference is Zhu, Miao and Peng [26]. These results are somewhat surprising from the viewpoint of nonparametric estimation. Note that, accordingly, the number of slices is similar to a tuning parameter such as, say, the bin width in a histogram estimator or, more generally, the bandwidth in a kernel estimator. We can regard a kernel estimator as a smoothed version of the slicing estimator with moving windows. However, as we know, to ensure $\sqrt{n}$ consistency of the kernel estimator, the bandwidth selection must be undertaken with care. Zhu and Fang [25] proved the asymptotic normality of the kernel estimator of SIR when the bandwidth is selected in the range $n^{-1/2}$ to $n^{-1/4}$, which means that in probability, each window must have $n^{\delta}$ points for some $\delta > 0$. Therefore, for SIR, Li's slicing estimation has the advantage that a less smoothed estimator is even less sensitive to the tuning parameter.

The problem of whether SAVE has similar properties to SIR is then of great interest. Empirical studies have examined this and there is a general feeling that SAVE may be more sensitive to the choice of the number of slices than SIR. Cook [5] mentioned that the number of slices plays the role of tuning parameter and thus SAVE may be affected by this choice. The empirical study of Zhu, Ohtaki and Li [28] was consistent with the sensitivity of SAVE to the selection of the number of slices, but no theoretical results have been produced to show why and how the number of slices affects the performance of SAVE.

In this paper, we present a systematic study of this problem and obtain the following results.

1. When $Y$ is discrete and takes a finite value, SAVE is able to achieve $\sqrt{n}$ consistency.
2. For continuous $Y$, the convergence of SAVE is almost completely different from that of SIR. Let $c$ denote the number of data points in each slice. When $c$ is a fixed constant, SAVE is not consistent. When $c \sim n^b$ with $b > 0$, although the estimator for SAVE is consistent, it cannot be $\sqrt{n}$ consistent.
3. A bias correction is proposed to allow the SAVE estimator to be $\sqrt{n}$ consistent. Since in practice, the discretized approximation is commonly used in the literature, we present asymptotic normality in a general setting.

Note that Cook and Ni ([8], Section 7) investigated the asymptotic behavior of the slicing estimator of the SAVE matrix and reported a result that is relevant to Theorem 2.3 in this paper. Another relevant paper is [12].

The rest of this paper is organized as follows. Section 2 contains an investigation into when the estimator is $\sqrt{n}$ consistent. Section 3 contains the bias correction and an approximation via discretization. Section 4 reports a simulation study and the performances of SIR, SAVE and the bias-corrected SAVE are considered. The proofs of the theorems are given in the Appendix.

**2. Asymptotic behavior of the slicing estimator.** As matrix operations are involved, we will write, unless stated otherwise, $\mathbf{A}\mathbf{A}^T = \mathbf{A}^2$, where $\mathbf{A}$ is a square matrix. We first describe the slicing estimator for the SAVE matrix $\mathbf{E}(I_p - \Sigma_{\boldsymbol{z}|y})^2$.

Suppose that $\{(\boldsymbol{z}_1, y_1), \dots, (\boldsymbol{z}_n, y_n)\}$ is a sample. Sort all of the data $(\boldsymbol{z}_i, y_i), i = 1, 2, \dots, n$, according to the ascending order of $y_i$. Define the order statistics $y_{(1)} \le y_{(2)} \le \cdots \le y_{(n)}$ and for every $1 \le i \le n$, let $\boldsymbol{z}_{(i)}$ be the concomitant of $y_{(i)}$. For any integer $c$, we group every $c$ data points and introduce a double subscript $(h, j)$, where $h$ refers to the slice number and $j$ refers to the order number of an observation in the given slice. Then

$$y_{(h,j)} = y_{(c(h-1)+j)}, \qquad \boldsymbol{z}_{(h,j)} = \boldsymbol{z}_{(c(h-1)+j)}, \qquad \bar{\boldsymbol{z}}_{(h)} = \frac{1}{c}\sum_{j=1}^{c} \boldsymbol{z}_{(h,j)}.$$

The number of data points in the last slice may be less than $c$, but the calculation is similar and the asymptotic results are still valid. Without loss of generality, suppose that we have $H$ slices and that $n = c \times H$. The sample version of the conditional variance of $\boldsymbol{z}$ given $y$ in each slice is

$$\hat{\Sigma}(h) = \frac{1}{(c-1)} \sum_{j=1}^{c} (\boldsymbol{z}_{(h,j)} - \bar{\boldsymbol{z}}_{(h)})^2. \tag{2.1}$$

The estimate of $\mathbf{E}((I_p - \Sigma_{\boldsymbol{z}|y})^2)$ is defined as

$$\frac{1}{H}\sum_{h=1}^{H}(I_p - \hat{\Sigma}(h))^2 = I_p - 2\frac{1}{H}\sum_{h=1}^{H}\hat{\Sigma}(h) + \frac{1}{H}\sum_{h=1}^{H}(\hat{\Sigma}(h))^2. \tag{2.2}$$

Note that the term $I_p - \frac{1}{H}\sum_{h=1}^{H}\hat{\Sigma}(h)$ is the same as the SIR estimator. Zhu and Ng [27] proved the $\sqrt{n}$ consistency of $I_p - \frac{1}{H}\sum_{h=1}^{H}\hat{\Sigma}(h)$ under certain regularity conditions. Hence, throughout the rest of the paper, we only investigate the asymptotic properties of $\Lambda_n = \frac{1}{H}\sum_{h=1}^{H}(\hat{\Sigma}(h))^2$, the results of the estimator of SAVE being presented as corollaries. Moreover, $\Lambda_n$ can be rewritten as

$$\begin{aligned}
\Lambda_n &= \frac{1}{H}\sum_{h=1}^{H}(\hat{\Sigma}(h))^2 \\
&= \frac{1}{H}\sum_{h=1}^{H}\left\{\frac{1}{(c-1)}\sum_{j=1}^{c}(\boldsymbol{z}_{(h,j)} - \bar{\boldsymbol{z}}_{(h)})^2\right\}^2 \\
&= \Bigg[\sum_{h=1}^{H}\sum_{l=2}^{c}\sum_{j=1}^{l-1}\sum_{v=2}^{c}\sum_{u=1}^{v-1}(\boldsymbol{z}_{(h,l)} - \boldsymbol{z}_{(h,j)})(\boldsymbol{z}_{(h,l)} - \boldsymbol{z}_{(h,j)})^T \\
&\qquad \times (\boldsymbol{z}_{(h,v)} - \boldsymbol{z}_{(h,u)})(\boldsymbol{z}_{(h,v)} - \boldsymbol{z}_{(h,u)})^T\Bigg][nc(c-1)^2]^{-1}.
\end{aligned}$$

For the sake of convenience, we here introduce some notation. For a symmetric $p\times p$ matrix $D = (d_{ij})$, $\mathrm{vech}\{D\} = (d^{(11)}, \ldots, d^{(p1)}, d^{(22)}, \ldots, d^{(p2)}, \ldots, d^{(pp)})^T$ is the $\frac{p(p+1)}{2}\times 1$ vector constructed from the elements of $D$.

We now define the total variation of order $r$ for a function. Let $\Pi_n(K)$ be the collection of $n$-point partitions $-K \le y_{(1)} \le \cdots \le y_{(n)} \le K$ of the closed interval $[-K, K]$, where $K > 0$ and $n \ge 1$. Any vector-valued or real-valued function $\mathbf{f}(y)$ is said to have a total variation of order $r$ if for any fixed $K > 0$,

$$\lim_{n\to\infty}\frac{1}{n^r}\sup_{\Pi_n(K)}\sum_{i=1}^{n}\|\mathbf{f}(y_{i+1}) - \mathbf{f}(y_i)\| = 0.$$

For any vector-valued or real-valued function $\mathbf{f}(y)$, if there are a nondecreasing real-valued function $M$ and a real number $K_0$ such that for any two points, say $y_1$ and $y_2$, both in $(-\infty, -K_0]$ or both in $[K_0, +\infty)$,

$$\|\mathbf{f}(y_1) - \mathbf{f}(y_2)\| \le |M(y_1) - M(y_2)|,$$

then we can say that the function $\mathbf{f}(y)$ is nonexpansive in the metric of $M$ on both sides of $K_0$.

2.1. *When is SAVE not $\sqrt{n}$ consistent?.* Let $\mathbf{m}(y) = \mathbf{E}(\boldsymbol{z}|Y = y)$. We can write $\boldsymbol{z} = \boldsymbol{\varepsilon} + \mathbf{m}(y)$, where $\mathbf{E}(\boldsymbol{\varepsilon}|Y) = 0$, and then $\Lambda = \mathbf{E}[(\Sigma_{\boldsymbol{z}|Y})^2] = \mathbf{E}[(\mathbf{E}(\boldsymbol{\varepsilon}\boldsymbol{\varepsilon}^T|Y))^2]$. The conditional expectation of $\boldsymbol{\varepsilon}$ given $y$ equals zero and more importantly, when $y_i$ are given, $\boldsymbol{\varepsilon}_i$ are independent, although they are not identically distributed (see [14] or [27]). Analogously to $\Lambda_n$, we denote

$$A_n = \left[\sum_{h=1}^{H}\sum_{l=2}^{c}\sum_{j=1}^{l-1}\sum_{v=2}^{c}\sum_{u=1}^{v-1}(\boldsymbol{\varepsilon}_{(h,l)} - \boldsymbol{\varepsilon}_{(h,j)})(\boldsymbol{\varepsilon}_{(h,l)} - \boldsymbol{\varepsilon}_{(h,j)})^T(\boldsymbol{\varepsilon}_{(h,v)} - \boldsymbol{\varepsilon}_{(h,u)}) \times (\boldsymbol{\varepsilon}_{(h,v)} - \boldsymbol{\varepsilon}_{(h,u)})^T\right][nc(c-1)^2]^{-1}.$$

Let $J_n = \Lambda_n - A_n$. To prove the convergence of $\Lambda_n$, we need to investigate $A_n$ and $J_n$.

THEOREM 2.1. *Assume the following four conditions:*

(1) *There is a nonnegative number $\alpha$ such that $\mathbf{E}(\|\boldsymbol{z}\|^{8+\alpha}) < \infty$.*

(2) *The inverse regression function $\mathbf{m}(y)$ has a total variation of order $r > 0$.*

(3) *$\mathbf{m}(y)$ is nonexpansive in the metric of $M(y)$ on both sides of a positive number $B_0$ such that*

$$M^{8+\alpha}(t)P(Y > t) \to 0 \qquad \textit{as } t \to \infty.$$

(4) *$c \sim n^b$ for $b \geq 0$.*

*Then $n^\beta J_n = o_p(1)$ for any $\beta$ such that $\beta + b + \max\{\frac{3}{8+\alpha} + r, \frac{4}{8+\alpha}\} \leq 1$.*

REMARK 2.1. We note that the conditions are similar to those that ensure the consistency of the estimator for SIR, except for the higher moments of $\boldsymbol{z}$ (see [27]). The $\sqrt{n}$ consistency of $J_n$ implies $\beta = 0.5$ and hence we must have $b = 1/2 - \max\{\frac{3}{8+\alpha} + r, \frac{4}{8+\alpha}\} \geq 0$. When $r$ is close to zero and all moments exist, $c$ can be selected to be arbitrarily close to $\sqrt{n}$.

THEOREM 2.2. *Assume the following conditions*:

(1) *There is a nonnegative number $\alpha$ such that $\mathbf{E}(\|\boldsymbol{z}\|^{\max\{8+\alpha,12\}}) < \infty$.*

(2) *Let $\mathbf{m}_1(y) = \mathbf{E}(\boldsymbol{\varepsilon}\boldsymbol{\varepsilon}^T|Y = y)$. $\mathbf{m}_1(y)$ has a total variation of order $r_1 > 0$.*

(3) *For a nondecreasing continuous function $M_1(\cdot)$, $\mathbf{m}_1(y)$ is nonexpansive in the metric of $M_1(y)$ on both sides of a positive number $B_0'$ such that*

$$M_1^{4+\alpha/2}(t)P(Y > t) \to 0 \qquad \textit{as } t \to \infty.$$

(4) *Let* $\mathbf{m}_2(y) = \mathbf{E}((\boldsymbol{\varepsilon}\boldsymbol{\varepsilon}^T)^2|y)$. *For a nondecreasing continuous function* $M_2(\cdot)$, $\mathbf{m}_2(y)$ *is nonexpansive in the metric of* $M_2(y)$ *on both sides of a positive number* $B_0''$ *such that*

$$M_2^{2+\alpha/4}(t)P(Y > t) \to 0 \qquad \text{as } t \to \infty.$$

(5) *There exists a positive* $\rho_1$ *such that*

$$\lim_{d\to\infty} \limsup_{n\to\infty} \mathbf{E}(|M_1^2(y_{(n)})|I(|M_1(y_{(n)})| > d)) = o(n^{-\rho_1}).$$

(6) *There exists a positive* $\rho_2$ *such that*

$$\lim_{d\to\infty} \limsup_{n\to\infty} \mathbf{E}(|M_2(y_{(n)})M_1^2(y_{(n)})|I(|M_2(y_{(n)})| > d)) = o(n^{-\rho_2}).$$

*Then*

$$\mathbf{E}(A_n) = \left(1 - \frac{(c-2)}{c(c-1)}\right)\Lambda + \frac{1}{c}\mathbf{E}[(\boldsymbol{\varepsilon}\boldsymbol{\varepsilon}^T)^2] + o(cn^{-1+\max\{r_1, \frac{2}{4+\alpha/2}, \rho_1\}}). \tag{2.3}$$

*On the further assumption that* $c \sim n^b$ *for* $b > 0$, *we have*

$$n^\beta(A_n - \Lambda) = o_p(1) \tag{2.4}$$

*for any* $\beta$ *such that* $\beta + b + \max\{r_1, \frac{2}{4+\alpha/2}, \rho_1\} \le 1$, $\beta < b$, *and* $2\beta + b + \max\{2r_1, \frac{2}{4+\alpha/2} + \frac{1}{2+\alpha/4}, \rho_2\} \le 2$.

REMARK 2.2. The first three conditions in Theorem 2.2 are similar to those in Theorem 2.1. Condition (2) is similar to the condition for the inverse regression function because we deal with the conditional second moment of $\boldsymbol{\varepsilon}$ when SAVE is applied. Condition (3) is slightly weaker than the existence of the $(4 + \alpha/2)$th moment of $M_1(\cdot)$ or, equivalently, the $(8 + \alpha)$th moment of $\boldsymbol{z}$, as is Condition (4). Note that Condition (5) is slightly stronger than $M_1^2(y_{(n)}) = o_p(n^{\rho_1})$ because we have to handle the moment convergence. It is well known that when the $y_i$ follow an exponential distribution, the maximum $y_{(n)}$ can be bounded by $(\log n)^c$ in probability for some $c \ge 1$ (see, e.g., [2], Chapter 1, page 10), and when the support of $y_i$ is bounded, $y_{(n)}$ is simply bounded by a constant. Note that for any transformation $h(\cdot)$ on $y$, $h(y)$ is independent of $\boldsymbol{z}$ when $B^T\boldsymbol{z}$ is given. Therefore, we could construct a transformation to allow the support of bounded $h(y)$ and consider the $(\boldsymbol{z}_i, h(y_i))$'s. However, in this paper we do not consider any transformations of $y$.

REMARK 2.3. From Theorems 2.1 and 2.2, we know that when $c$ is a fixed constant, $J_n = o_p(1)$, but the mean of $A_n$ is not asymptotically equal to $\Lambda$. From the proof of Theorem 2.2, we can easily see that $A_n$ does not converge in probability to $\Lambda$ and therefore $\Lambda_n = J_n + A_n$ cannot converge to $\Lambda$. When $c$ tends to infinity at a rate slower than $n^{1/2}$ in Theorems 2.1 and 2.2, the convergence rate of $\Lambda_n$ to $\Lambda$ is slower than $1/c$ and therefore $\sqrt{n}$ consistency does not hold. This property is completely different from that of SIR because within this range of $c$, the slicing estimator of SIR is $\sqrt{n}$ consistent (see [27]). The second and third terms in $\mathbf{E}(A_n)$ provide two bounds, when $r_1 = 0$, $\alpha = \infty$ with the multiplication of $\sqrt{n}$ by $\mathbf{E}(A_n)$, $\sqrt{n}/c$ and $c/\sqrt{n}$, that are reciprocal one to another. Although the third term is an upper bound, it is tight, to a certain extent. An example is provided by the case where $y$ is uniformly distributed on [0,1], $y_{(i)} = i/n$. With large probability so the third term can achieve the rate $cn^{-1}$, which means that in general cases, if no extra conditions are imposed, it is impossible for the expectation of $A_n$ to converge to $\Lambda$. This can be seen from the proof of the theorem. This is worthy of a detailed investigation and relates to the question of whether the slicing estimator of SAVE is $\sqrt{n}$ consistent. In the following subsection, we undertake a detailed study of this issue.

When the mean and covariance of $\boldsymbol{x}$ are unknown, the $\hat{\boldsymbol{z}}_i = \Sigma_{\boldsymbol{x}}^{-1/2}(\boldsymbol{x}_i - \bar{\boldsymbol{x}})$ are used to estimate the matrix $\mathbf{E}(I_p - \Sigma_{\boldsymbol{z}|Y})^2$. Let $\hat{\Sigma}_{\hat{\boldsymbol{z}}}(h)$ be the sample covariance of the $\hat{\boldsymbol{z}}_i$'s in each slice for $h = 1, \ldots, H$. Note that this matrix is location-invariant. We can assume, with no loss of generality, that the sample mean $\bar{\boldsymbol{x}} = 0$. Clearly, $\hat{\Sigma}_{\hat{\boldsymbol{z}}}(h) = \hat{\Sigma}_{\boldsymbol{x}}^{-1/2}\Sigma_{\boldsymbol{x}}{}^{1/2}\hat{\Sigma}(h)\Sigma_{\boldsymbol{x}}{}^{1/2}\hat{\Sigma}_{\boldsymbol{x}}^{-1/2}$. To study the asymptotic behavior of the estimator when $\Sigma_{\boldsymbol{x}}$ is replaced by $\hat{\Sigma}_{\boldsymbol{x}}$, we first consider the following property. Let $R = (\hat{\Sigma}_{\boldsymbol{x}} - \Sigma_{\boldsymbol{x}})\Sigma_{\boldsymbol{x}}^{-1}$. By some elementary calculation and the well-known fact that $\hat{\Sigma}_{\boldsymbol{x}} - \Sigma_{\boldsymbol{x}} = O_p(1/\sqrt{n})$, we have

$$\begin{aligned}
\hat{\Sigma}_{\boldsymbol{x}}^{-1/2}\Sigma_{\boldsymbol{x}}{}^{1/2} &= I_p - (\hat{\Sigma}_{\boldsymbol{x}} - \Sigma_{\boldsymbol{x}})\Sigma_{\boldsymbol{x}}^{-1}[(I_p + R)^{-1}((I_p + R)^{-1/2} + I_p)^{-1}] \\
&= I_p - \frac{1}{2}(\hat{\Sigma}_{\boldsymbol{x}} - \Sigma_{\boldsymbol{x}})\Sigma_{\boldsymbol{x}}^{-1} + o_p(1/\sqrt{n})
\end{aligned} \tag{2.5}$$

and similarly

$$\Sigma_{\boldsymbol{x}}^{1/2}\hat{\Sigma}_{\boldsymbol{x}}^{-1/2} = I_p - \frac{1}{2}\Sigma_{\boldsymbol{x}}^{-1}(\hat{\Sigma}_{\boldsymbol{x}} - \Sigma_{\boldsymbol{x}}) + o_p(1/\sqrt{n}). \tag{2.6}$$

Consequently, for each $h = 1, \ldots, H$,

$$\begin{aligned}
&\hat{\Sigma}_{\boldsymbol{x}}^{-1/2}\Sigma_{\boldsymbol{x}}{}^{1/2}\hat{\Sigma}(h)\Sigma_{\boldsymbol{x}}{}^{1/2}\hat{\Sigma}_{\boldsymbol{x}}^{-1/2} \\
&\quad = \hat{\Sigma}(h) - \frac{1}{2}(\hat{\Sigma}_{\boldsymbol{x}} - \Sigma_{\boldsymbol{x}})\Sigma_{\boldsymbol{x}}^{-1}\hat{\Sigma}(h) - \frac{1}{2}\hat{\Sigma}(h)\Sigma_{\boldsymbol{x}}^{-1}(\hat{\Sigma}_{\boldsymbol{x}} - \Sigma_{\boldsymbol{x}}) + o_p(1/\sqrt{n})
\end{aligned} \tag{2.7}$$

and then

$$\frac{1}{H}\sum_{h=1}^{H}(I_p - \hat{\Sigma}_{\boldsymbol{x}}^{-1/2}\Sigma_{\boldsymbol{x}}{}^{1/2}\hat{\Sigma}(h)\Sigma_{\boldsymbol{x}}{}^{1/2}\hat{\Sigma}_{\boldsymbol{x}}^{-1/2})^2$$

$$
\begin{aligned}
&= \frac{1}{H}\sum_{h=1}^{H}(I_p-\hat{\Sigma}(h))^2 \\
&\quad + \frac{1}{2H}\sum_{h=1}^{H}[(\hat{\Sigma}_{\boldsymbol{x}}-\Sigma_{\boldsymbol{x}})\Sigma_{\boldsymbol{x}}^{-1}\hat{\Sigma}(h)+\hat{\Sigma}(h)\Sigma_{\boldsymbol{x}}^{-1}(\hat{\Sigma}_{\boldsymbol{x}}-\Sigma_{\boldsymbol{x}})](I_p-\hat{\Sigma}(h)) \\
&\quad + \frac{1}{2H}\sum_{h=1}^{H}(I_p-\hat{\Sigma}(h))[(\hat{\Sigma}_{\boldsymbol{x}}-\Sigma_{\boldsymbol{x}})\Sigma_{\boldsymbol{x}}^{-1}\hat{\Sigma}(h)+\hat{\Sigma}(h)\Sigma_{\boldsymbol{x}}^{-1}(\hat{\Sigma}_{\boldsymbol{x}}-\Sigma_{\boldsymbol{x}})] \\
&\quad + o_p(1/\sqrt{n}\,) \\
&=: \frac{1}{H}\sum_{h=1}^{H}(I_p-\hat{\Sigma}(h))^2 + I_n + o_p(1/\sqrt{n}).
\end{aligned} \tag{2.8}
$$

We now deal with $I_n$. Write $(\hat{\Sigma}_{\boldsymbol{x}}-\Sigma_{\boldsymbol{x}})\Sigma_{\boldsymbol{x}}^{-1}=A_n=(a_{n,ij})$, $\hat{\Sigma}(h)=B_n(h)=(b_{n,ij}(h))$ and $(I_p-\hat{\Sigma}(h))=C_n(h)=(c_{n,ij}(h))$. $\sqrt{n}I_n$ can be written as

$$
\sqrt{n}I_n=\frac{\sqrt{n}}{2H}\sum_{h=1}^{H}[(A_nB_n(h)+B_n(h)A_n^T)C_n(h)+C_n(h)(A_nB_n(h)+B_n(h)A_n^T)]
$$

and its elements have the formula

$$
\begin{aligned}
\sqrt{n}I_{nil} &= \sum_{k=1}^{p}\sum_{j=1}^{p}\sqrt{n}a_{nlk}\frac{1}{2H}\sum_{h=1}^{H}[b_{njk}(h)c_{nkl}(h)+c_{nij}(h)b_{njk}(h)] \\
&\quad + \sum_{k=1}^{p}\sum_{j=1}^{p}\sqrt{n}a_{nkj}\frac{1}{2H}\sum_{h=1}^{H}[b_{nji}(h)c_{nlk}(h)+b_{nkl}(h)c_{nji}(h)] \\
&=: \sum_{k=1}^{p}\sum_{j=1}^{p}\sqrt{n}a_{nlk}D_{nijkl}.
\end{aligned} \tag{2.9}
$$

From the proofs of Theorems 2.1 and 2.2 in the Appendix, $D_{n\,ijkl}$ converges in probability to a constant $\tilde{D}_{ijkl}$. The well-known result of sample covariance yields the asymptotic normality of all $\sqrt{n}a_{nil}$. Thus, $\sqrt{n}I_{nil}$ converges in distribution to $N(0,V_{il})$, where $V_{il}=\lim_{n\to\infty}\operatorname{var}(\sum_{k=1}^{p}\sum_{j=1}^{p}\sqrt{n}a_{nlk}\tilde{D}_{ijkl})$. This means that $I_{nil}=O_p(1/\sqrt{n})$ and we have the following result.

COROLLARY 2.1. *Under the conditions of Theorems 2.1 and 2.2, the results of these two theorems continue to hold when the mean and covariance of $\boldsymbol{x}$ are unknown and the $\hat{\boldsymbol{z}}_i=\Sigma_{\boldsymbol{x}}^{-1/2}(\boldsymbol{x}_i-\bar{\boldsymbol{x}})$ are used to estimate the matrix* $\mathbf{E}(I_p-\Sigma_{\boldsymbol{z}|Y})^2$.

This corollary holds because the convergence rate of $I_n$ is faster than the convergence rate of $\Lambda_n$ and thus the results of Theorems 2.1 and 2.2 do not change.

2.2. *When is SAVE $\sqrt{n}$ consistent?.* The following theorem asserts the asymptotic normality of the estimator in a special case in which the response is discrete and takes a finite value. For any value $l$, define $E_1(l) = \mathbf{E}(\boldsymbol{z}|Y=l)$ and

$$V(Y,\boldsymbol{z}) = \sum_{l=1}^{d}[-2((\boldsymbol{z}_j^2 - 2\boldsymbol{z}_j E_1(l))I(y_j=l) - \mathbf{E}((\boldsymbol{z}^2 - 2\boldsymbol{z}E_1(l))I(Y=l))) \\ \times (I_p - \mathrm{Cov}(\boldsymbol{z}|Y=l)) + (I(y_j=l)-p_l)\times(I_p - \mathrm{Cov}(\boldsymbol{z}|Y=l))^2].$$

THEOREM 2.3. *Assume that the response $Y$ takes $d$ values and, without loss of generality, assume that $Y = 1,2,\ldots,d$ and $P(Y=l) = p_l > 0$ for $l = 1,\ldots,d$. Additionally, assume that $\mathbf{E}\|\boldsymbol{z}\|^8 < \infty$. Then when $H = d$,*

$$\sqrt{n}\,\mathrm{vech}\left(\frac{1}{H}\sum_{h=1}^{H}(I_p - \hat{\Sigma}(h))^2 - \mathbf{E}(I_p - \Sigma_{\boldsymbol{z}|Y})^2\right) \Rightarrow N(0, \mathrm{Cov}(\mathrm{vech}\{V(Y,\boldsymbol{z})\}).$$

When the $\hat{\boldsymbol{z}}_j$ are used to estimate the SAVE matrix, the term $\sqrt{n}I_n$ affects the limiting variance. Note that

$$(\hat{\Sigma}_{\boldsymbol{x}} - \Sigma_{\boldsymbol{x}})\Sigma_{\boldsymbol{x}}^{-1} = \frac{1}{n}\sum_{j=1}^{n}[(\boldsymbol{x}_j - E(\boldsymbol{x}))^2 - \Sigma_{\boldsymbol{x}}]\Sigma_{\boldsymbol{x}}^{-1} + o_p(1/\sqrt{n}\,) \\ =: \frac{1}{n}\sum_{m=1}^{n}(e_{mlk})_{1\le k,\ l\le p} + o_p(1/\sqrt{n}\,). \tag{2.10}$$

The leading term is a sum of i.i.d. random variables, which implies that $a_{nlk}$ is asymptotically a sum of i.i.d. random variables. Then from (2.9),

$$\sqrt{n}(I_{nil})_{1\le i,\ l\le p} = \frac{1}{\sqrt{n}}\sum_{m=1}^{n}\left(\sum_{k=1}^{p}\sum_{j=1}^{p}e_{mlk}D_{nijkl}\right)_{1\le i,\,l\le p} + o_p(1) \\ =: \frac{1}{\sqrt{n}}\sum_{m=1}^{n}\mathbf{E}_m + o_p(1). \tag{2.11}$$

COROLLARY 2.2. *Under the conditions of Theorem 2.3,*

$$\sqrt{n}\,\mathrm{vech}\left(\frac{1}{H}\sum_{h=1}^{H}(I_p - \hat{\Sigma}_{\hat{\boldsymbol{z}}}(h))^2 - \mathbf{E}(I_p - \Sigma_{\boldsymbol{z}|Y})^2\right) \\ \Rightarrow N(0, \mathrm{Cov}(\mathrm{vech}\{V(Y,\boldsymbol{z}) + \mathbf{E}_1\}).$$

## 3. The approximation and bias correction.

3.1. *The approximation.* Note that when $Y$ is a discrete random variable, SAVE needs only very mild conditions to achieve asymptotic normality. In this case, $H$ is a fixed number that does not depend on $n$. In applications, $H$ is often a fixed number, which means that approximation via discretization is used in practice. It would be worthwhile to conduct a theoretical investigation to ascertain the rationale of the approximation.

Let $S_h = (q_{h-1}, q_h]$ for $h = 1, \dots, H$, $q_0 = -\infty$, $q_H = \infty$ and $p_h = P(Y \in S_h)$. Recall that the construction of the slicing estimator is based on a weighted sum of the sample covariance matrices of the associated $\boldsymbol{z}_i$'s with $y_i$'s in all slices $S_h$, $h = 1, \dots, H$. These sample covariance matrices are the estimators of the $\mathbf{E}(\mathrm{Cov}(\boldsymbol{z}|Y \in S_h))$'s. Note that these matrices can be written as

$$\Sigma(h) := \frac{\mathbf{E}((\boldsymbol{z} - \frac{\mathbf{E}(\boldsymbol{z} I(Y \in S_h))}{p_h})^2 I(Y \in S_h))}{p_h},$$

where $I(\cdot)$ is the indicator function. The estimator of $p_h$ is equal to $1/H$ when $q_h$ is replaced by the empirical quantile $\hat{q}_h$. The slicing estimator can be rewritten as $I_p - \frac{2}{H}\sum_{h=1}^H \hat{\Sigma}(h) + \frac{1}{H}\sum_{h=1}^H \hat{\Sigma}^2(h)$
with

$$\begin{aligned}\hat{\Sigma}(h) &= \frac{1}{c}\sum_{j=1}^{c}(\boldsymbol{z}_{(h,j)} - \bar{\boldsymbol{z}}_{(h)})^2 \\ &= \frac{1}{n\hat{p}_h}\sum_{j=1}^{n}\left(\boldsymbol{z}_j - \frac{1}{n\hat{p}_h}\sum_{j=1}^{n}\boldsymbol{z}_j I(y_j \in \hat{S}_h)\right)^2 I(y_j \in \hat{S}_h).\end{aligned} \tag{3.1}$$

That is, the slicing estimator estimates $\Lambda(H) = \sum_{h=1}^H (I_p - \Sigma(h))^2 p_h$. In the case in which $Y$ is continuous and $H$ is large, we have

$$\begin{aligned}\Lambda(H) &\cong \sum_{h=1}^{H} \mathbf{E}[(I_p - \mathrm{Cov}(\boldsymbol{z}|Y))^2 I(Y \in S_h)] \\ &= \mathbf{E}(I_p - \mathrm{Cov}(\boldsymbol{z}|Y))^2,\end{aligned}$$

where $\cong$ stands for approximate equality. Clearly, under some regularity conditions, $\Lambda(H)$ can converge to $\mathbf{E}((I_p - \mathrm{Cov}(\boldsymbol{z}|Y))^2)$ as $H \to \infty$.

As with Theorem 2.3, we have the following result. Define, for every $h$, $E_1(h) = \mathbf{E}(\boldsymbol{z}|Y \in S_h)$ and take $f(q_j)$ as being the value of the density of $Y$ at $q_j$.

THEOREM 3.1. *Let $\hat{q}_h = y_{(ch)}$, $h = 1, \dots, H-1$, be the empirical $(h/H)$th quantiles, with $\hat{q}_0 = 0$ and $\hat{q}_H = \infty$. Assume the following:*

(1) $\mathbf{E}\|\boldsymbol{z}\|^8 < \infty$.

(2) *If we write* $\mathbf{E}(F(Y, \boldsymbol{z}, a, b)) := \mathbf{E}(\boldsymbol{z}^2(I(Y \in (a, b]) - I(Y \in S_h))$, *then* $\mathbf{E}(F(Y, \boldsymbol{z}, a, b))$ *is differentiable with respect to* $a$ *and* $b$ *and its first derivative is bounded by a constant* $C_1$.

(3) *If we write* $\mathbf{E}(G(Y, \boldsymbol{z}, a, b)) := \mathbf{E}(\boldsymbol{z}(I(Y \in (a, b]) - I(Y \in S_h)))$, *then* $\mathbf{E}(G(Y, \boldsymbol{z}, a, b))$ *is differentiable with respect to* $a$ *and* $b$.

(4) *The density function* $f(y)$ *of* $Y$ *is bounded away from zero at all quantiles* $q_h$, $h = 1, \ldots, H-1$.

*When* $\Lambda_n$ *is constructed with the slices* $\hat{S}_h = (\hat{q}_{h-1}, \hat{q}_h]$, $h = 1, \ldots, H$, *as* $n \to \infty$,

$$\sqrt{n}\,\mathrm{vech}\left(\frac{1}{H}\sum_{h=1}^{H}(I_p - \hat{\Sigma}(h))^2 - \mathbf{E}(I_p - \Sigma_{\boldsymbol{z}|Y})^2\right)$$

*is asymptotically normal with zero mean and variance* $\mathrm{Cov}(\mathrm{vech}\{L(Y, \boldsymbol{z})\})$. *When the* $\hat{\boldsymbol{z}}_i$ *are used to construct the estimator, the limiting variance is* $\mathrm{Cov}(\mathrm{vech}\{L(Y, \boldsymbol{z}) + \mathbf{E}_1\})$, *where*

$$\begin{aligned} L(Y, \boldsymbol{z}) = \Biggl\{ &-2\sum_{h=1}^{H}((\boldsymbol{z}^2 - 2\boldsymbol{z}E_1(h))I(Y \in S_h) - \mathbf{E}((\boldsymbol{z}^2 - 2\boldsymbol{z}E_1(h))I(Y \in S_h))) \\ &- 2\sum_{h=1}^{H}\left(\frac{-I(Y \le q_{h-1}) + \frac{h-1}{H}}{f(q_{h-1})}, -\frac{I(Y \le q_h) + \frac{h}{H}}{f(q_h)}\right) \\ &\times (\tilde{F}'(q_{h-1}, q_h) - 2\tilde{G}'(q_{h-1}, q_h)E_1(h))\Biggr\} \times (I_p - \Sigma(h)) \end{aligned}$$

*and* $\mathbf{E}_1$ *is defined as in* (2.11).

Remark 3.1. Conditions (2)–(4) are assumed in order to ensure some degree of smoothness of the relevant functions, and thus the conditions are fairly mild.

3.2. *Bias correction.* In terms of examining the expectation of $A_n$, we can see that the major bias is the term $\frac{1}{c-1}\mathbf{E}(\boldsymbol{\varepsilon}\boldsymbol{\varepsilon}^T)^2$. If we can eliminate the impact of this term, then asymptotic normality may be possible. In this subsection, we suggest a bias correction, the idea of which is simple. We first obtain an estimator of this term and then subtract it from the estimator of $\Lambda_n$, which motivates the bias correction as follows.

As before, we divide the range of $Y$ into $H$ slices. According to the result of Theorem 2.2, the estimator of $V =: \mathbf{E}(\boldsymbol{\varepsilon}\boldsymbol{\varepsilon}^T)^2$ is defined as

$$V_n = \frac{1}{Hc}\sum_{h=1}^{H}\sum_{j=1}^{c}((\boldsymbol{z}_{(h,j)} - \bar{\boldsymbol{z}}_{(h)})(\boldsymbol{z}_{(h,j)} - \bar{\boldsymbol{z}}_{(h)})^T)^2.$$

The corrected estimator of $\Lambda$ is

$$\tilde{\Lambda}_n = \frac{c(c-1)}{(c-1)^2+1}\Lambda_n - \frac{c-1}{(c-1)^2+1}V_n.$$

THEOREM 3.2. *Assume that conditions* (2)–(3) *of Theorem* 2.1 *and conditions* (1)–(6) *of Theorem* 2.2 *are satisfied. Let* $c \sim n^b$, *where* $b$ *is a positive number that satisfies the following three inequalities*:

(a) $b > \frac{1}{4}$;
(b) $b \le 0.5 - \max\{\rho_1, r_1, \frac{2}{4+\alpha/2}, \frac{3}{8+\alpha} + r, \frac{4}{8+\alpha}\}$;
(c) $b \le 1 - \max\{2r_1, \frac{2}{4+\alpha/2} + \frac{1}{2+\alpha/4}, \rho_2\}$.

*Then* $\operatorname{vech}\frac{\sqrt{n}}{c}(V_n - V) = o_p(1)$ *and therefore* $\sqrt{n}\operatorname{vech}(\tilde{\Lambda}_n - \Lambda) = O_p(1)$. *The results continue to hold when the* $\hat{\boldsymbol{z}}_i$*'s are used to construct the estimators.*

Similarly to (2.9), the term that relates to $\hat{\Sigma}_{\boldsymbol{x}} - \Sigma_{\boldsymbol{x}} = O_p(1/\sqrt{n})$ and the $V_n$ that is based on the $\hat{\boldsymbol{z}}_i$'s differs by a term that is $O_p(1/\sqrt{n})$ from the $V_n$ that is based on the $\boldsymbol{z}_i$'s. Thus, the estimators that are based on the $\hat{\boldsymbol{z}}_i$'s have the same asymptotic behavior as that of the $V_n$ that are based on the $\boldsymbol{z}_i$'s.

To show the $\sqrt{n}$ consistency of the estimated CDR subspace, we define a bias-corrected estimator for the matrix $\mathbf{E}(I_p - \Sigma_{\boldsymbol{z}|y})^2$ by

$$\mathrm{CSAVE}_n := I_p - \frac{2}{H}\sum_{h=1}^{H}\hat{\Sigma}(h) + \tilde{\Lambda}_n.$$

The eigenvectors that are associated with the largest $k$ eigenvalues of $\mathrm{CSAVE}_n$ are used to form a basis of the estimated CDR space. following result asserts the asymptotic normality of the corrected estimator.

COROLLARY 3.1. *Under the conditions of Theorem* 3.2,

$$\sqrt{n}\operatorname{vech}(\mathrm{CSAVE}_n - \mathbf{E}((I_p - \Sigma_{\boldsymbol{z}|Y})^2))$$

*is asymptotically multinormal with zero mean and finite variance* $(\Delta_1 + \Delta_2)$, *where* $\Delta_1$ *and* $\Delta_2$ *are defined in* (A.17) *and* (A.19), *respectively. When the* $\hat{\boldsymbol{z}}_i$ *are used to construct* $\mathrm{CSAVE}_n$, *the limiting variance is* $(\Delta_1 + \Delta_2 + \mathbf{E}_1)$, *where* $\mathbf{E}_1$ *is the random matrix that is defined in* (2.11).

3.3. *The consistency of estimated eigenvalues and eigenvectors.* As the CDR space is estimated by the space that is spanned by the eigenvectors that are associated with the nonzero eigenvalues of the estimated SAVE matrix, we present the convergence of the estimated eigenvalues and eigenvectors. Because the convergence is the direct extension of the results of Zhu and Ng [27] or Zhu and Fang [25], we do not give the details of the proof in this paper.

From the theorems and corollary in this section, we can derive the asymptotic normality of the eigenvalues and the corresponding eigenvectors by using perturbation theory. The following result is parallel to the result for SIR obtained by Zhu and Fang [25] and Zhu and Ng [27]. The proof is also almost identical to that for the SIR matrix estimator. We omit the details of the proof in this article.

Let $\lambda_1(A) \geq \lambda_2(A) \geq \cdots \geq \lambda_p(A) \geq 0$ and $b_i(A) = (b_{1i}(A), \ldots, b_{pi}(A))^T$, $i = 1, \ldots, p$, denote the eigenvalues and their corresponding eigenvectors for a $p \times p$ matrix $A$. Let $\tilde{\Lambda} = \mathbf{E}(I_p - \Sigma_{\boldsymbol{z}|y})^2$ and $\bar{\Lambda}_n$ be the estimator that is defined in the theorems and corollary of Section 3.

THEOREM 3.3. *In addition to the conditions of the respective theorems in this section, assume that the nonzero $\lambda_l(\bar{\Lambda})$'s are distinct. Then for each nonzero eigenvalue $\lambda_i(\bar{\Lambda})$ and the corresponding eigenvector $b_i(\Lambda)$, we have*

$$\begin{aligned}\sqrt{n}(\lambda_i(\bar{\Lambda}_n) - \lambda_i(\bar{\Lambda}))& \\ &= \sqrt{n} b_i(\bar{\Lambda})^T (\bar{\Lambda}_n - \bar{\Lambda}) b_i(\bar{\Lambda}) + o_p(\sqrt{n}\|\bar{\Lambda}_n - \bar{\Lambda}\|) \\ &= b_i(\bar{\Lambda})^T W b_i(\bar{\Lambda}),\end{aligned} \tag{3.2}$$

*where $W$ is the limit matrix of $\sqrt{n}(\mathrm{CSAVE}_n - \mathbf{E}((I_p - \Sigma_{\boldsymbol{z}|Y})^2))$ that is studied in Corollary 3.1, and as $n \to \infty$,*

$$\begin{aligned}\sqrt{n}(b_i(\bar{\Lambda}_n) - b_i(\bar{\Lambda}))& \\ &= \sqrt{n} \sum_{l=1, l\neq i}^{p} \frac{b_i(\bar{\Lambda}) b_i(\bar{\Lambda})^T (\bar{\Lambda}_n - \bar{\Lambda}) b_i(\bar{\Lambda})}{\lambda_j(\bar{\Lambda}) - \lambda_l(\bar{\Lambda})} + o_p(\sqrt{n}\|\bar{\Lambda}_n - \bar{\Lambda}\|) \\ &= \sum_{l=1, l\neq i}^{p} \frac{b_i(\bar{\Lambda}) b_i(\bar{\Lambda})^T W b_i(\bar{\Lambda})}{\lambda_j(\bar{\Lambda}) - \lambda_l(\bar{\Lambda})},\end{aligned} \tag{3.3}$$

*where $\|\bar{\Lambda}_n - \bar{\Lambda}\| = \sum_{1 \leq i,j \leq p} |a_{ij}|$.*

**4. Simulation study and applications.** In this section, a simulation study is carried out to provide evidence for the efficiency of SIR, SAVE and the bias-corrected SAVE in practice. Following Li [16], the correlation coefficient between two spaces is taken to be the measure of the distance between the estimated CDR space and the true CDR space $S_{y|\boldsymbol{z}}$. For any eigenvector $\hat{\beta}_1$ that is associated with one of the largest $k$ eigenvalues obtained by the estimate, the squared multiple correlation coefficient $R^2(\hat{\beta}_1)$ between $\hat{\beta}_1^T \boldsymbol{z}$ and the ideally reduced variables $\beta_1^T \boldsymbol{z}, \ldots, \beta_k^T \boldsymbol{z}$ of $S_{y|\boldsymbol{z}}$ is employed to measure the distance between $\hat{\beta}_1$ and the space $S_{y|\boldsymbol{z}}$. That is,

$$R^2(\hat{\beta}_1) = \max_{\beta \in S_{y|\boldsymbol{z}}} \frac{(\hat{\beta}_1^T \Sigma_{\boldsymbol{z}} \beta)^2}{\hat{\beta}_1^T \Sigma_{\boldsymbol{z}} \hat{\beta}_1 \cdot \beta^T \Sigma_{\boldsymbol{z}} \beta}.$$

As $\boldsymbol{z}$ is a standardized variable, $R^2(\hat{\beta}_1)$ actually has the simpler formula

$$R^2(\hat{\beta}_1) = \max_{\beta \in S_{y|\boldsymbol{z}}} (\hat{\beta}_1^T \beta)^2.$$

When the estimated CDR space has dimension $k$, for a collection of the $k$ eigenvectors $\hat{\beta}_i$, $i = 1, \ldots, k$, that are associated with the $k$ largest eigenvalues, we use the squared trace correlation [the average of the squared canonical correlation coefficients between $\hat{\beta}_1^T \boldsymbol{z}, \ldots, \hat{\beta}_k^T \boldsymbol{z}$ and $\beta_1^T \boldsymbol{z}, \ldots, \beta_k^T \boldsymbol{z}$ as denoted by $R^2(\hat{\mathcal{B}})$] as our criterion (see also [13]), where $\tilde{\mathcal{B}}$ is the space that is spanned by $\{\hat{\beta}_1, \ldots, \hat{\beta}_k\}$.

We consider the cases where $k = 1$ and $n = 200$ and 480 and choose the following five models:

*Model* 1: $y = (\beta^T \boldsymbol{z})^3 + \varepsilon$.
*Model* 2: $y = (\beta^T \boldsymbol{z})^2 + \varepsilon$.
*Model* 3: $y = \beta^T \boldsymbol{z} \times \varepsilon$.
*Model* 4: $y = (\beta^T \boldsymbol{z})^3 + (\beta^T \boldsymbol{z}) \times \varepsilon$.
*Model* 5: $y = \cos(\beta^T \boldsymbol{z}) + \varepsilon$.

In these models, the covariate $\boldsymbol{z}$ and the error $\varepsilon$ are independent and respectively follow the normal distributions $N(0, I_{10})$ and $N(0, 1)$, where $I_{10}$ is the $10 \times 10$ identity matrix. In performing the simulation, we set $\beta = (1, 0, \ldots, 0)$.

We select models 1 to 5 based on the following considerations. Model 1 favors SIR rather than SAVE because the regression functions are strictly increasing. A similar investigation was undertaken in [28]. Model 2 favors SAVE rather than SIR because the inverse regression function is a zero function and then $\dim(S_{E(\boldsymbol{z}|y)}) = 0$ where $\dim(S)$ stands for the dimension of the space $S$. Model 3 deals with the variance function. Model 4 is constructed to be a combination of Model 1 and Model 3, as we are curious about the performance of SIR and SAVE in relation to the mean function and the variance function. We also include Model 5, which involves a periodic function.

The results are reported in Figure 1 and Table 1. When $n = 200$, a simulation was conducted with $H = 2$, 5, 10, 20 and 50, but we only report the results with $H = 10$ for illustration because for practical use, $H = 10$ is a good choice for this sample size (see relevant references such as [5, 16, 28]). The sensitivity to the slice selection will be discussed in terms of the results that are reported in Table 1 with $n = 480$. The boxplots in Figure 1 show the distribution of $R^2$ for a total of 200 Monte Carlo samples and show how the bias correction works with a fairly small sample size. From Figure 1, it is clear that CSAVE works well and is robust against the models that we employ.

Table 1 displays the numerical results for $n = 480$. The median of $R^2$ from a total of 200 Monte Carlo samples is presented so that we can compare the efficiency of the methods. To check the impact of the number of slices $H$, the values 2, 6, 24 and 96 are considered.

As expected, SIR is insensitive to $c$, but sensitive to the model and does not work well when the regression function is even or the CDR space is related to the error term.

The performance of SAVE is strongly affected by the choice of $c$, but when $H$ is properly chosen, SAVE works very well. However, the range of $c$ that results in a good performance from SAVE is fairly narrow. From the simulation results, we can see that when $H = 96$, that is, when $c = 5$, SAVE does not perform well. This is consistent with the theoretical conclusions in Section 2. The simulations show that choosing a relatively small $H$ favors SAVE, but that CSAVE still outperforms SAVE. Specifically, for $H = 2$, 6,

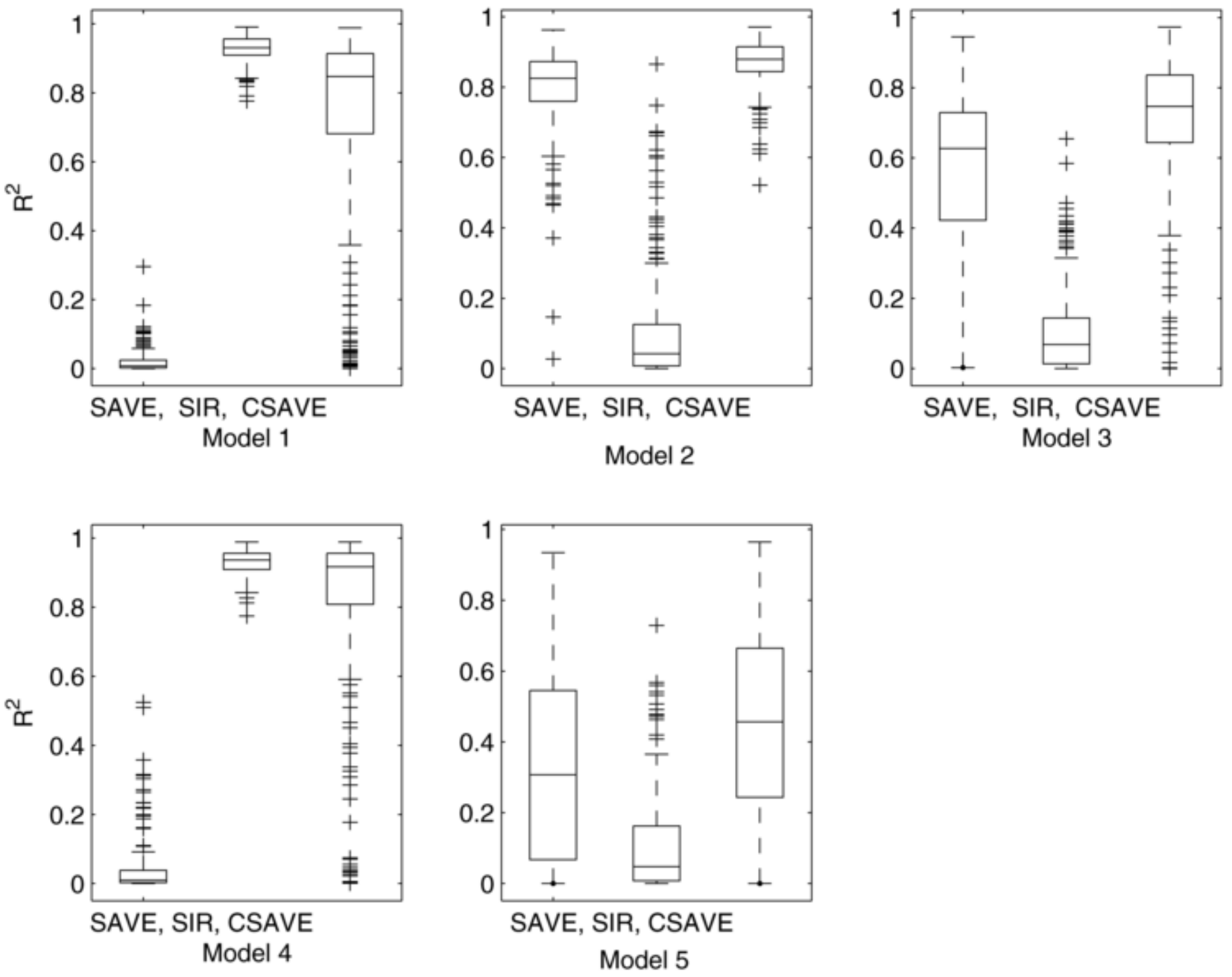

Fig. 1. *Boxplots of the distribution of* 200 *replicates of the* $R^2$ *values for models* 1–5 *when* $H = 10$ *and* $n = 200$. *The boxplots are, from left to right, for SAVE, SIR and CSAVE.*

TABLE 1
*The empirical median of the $R^2$ with $n = 480$*

| | $R^2(\hat{\beta})$ | | | |
|---|---|---|---|---|
| | $H = 2$ | $H = 6$ | $H = 24$ | $H = 96$ |
| *Model* 1 | | | | |
| SAVE | 0.7521 | 0.9599 | 0.0099 | 0.0009 |
| SIR | 0.9442 | 0.9681 | 0.9714 | 0.9586 |
| CSAVE | 0.8023 | 0.9687 | 0.9539 | 0.0122 |
| *Model* 2 | | | | |
| SAVE | 0.9539 | 0.9523 | 0.9187 | 0.7225 |
| SIR | 0.0460 | 0.0386 | 0.0443 | 0.0435 |
| CSAVE | 0.9575 | 0.9584 | 0.9317 | 0.8487 |
| *Model* 3 | | | | |
| SAVE | 0.0724 | 0.9201 | 0.8517 | 0.3547 |
| SIR | 0.0586 | 0.0545 | 0.0564 | 0.0448 |
| CSAVE | 0.0654 | 0.9336 | 0.8854 | 0.6393 |
| *Model* 4 | | | | |
| SAVE | 0.0741 | 0.9055 | 0.8665 | 0.3059 |
| SIR | 0.8656 | 0.8952 | 0.8825 | 0.7263 |
| CSAVE | 0.1066 | 0.9277 | 0.9024 | 0.7097 |
| *Model* 5 | | | | |
| SAVE | 0.8750 | 0.8657 | 0.6741 | 0.1249 |
| SIR | 0.0581 | 0.0484 | 0.0558 | 0.0625 |
| CSAVE | 0.8851 | 0.8966 | 0.7639 | 0.2517 |

24 and 96, the $R^2$ of CSAVE is larger than that of SAVE, especially when $H$ is large. Although the performance of CSAVE is also influenced by the choice of $c$, the range of $c$ that makes CSAVE work well is larger than that which makes SAVE work well. As, to some extent, CSAVE removes uncertainties about which $c$ should be used in practice, we recommend this method. Based on the limited simulations, $H = n/20$ is recommended for practical use.

## APPENDIX

As the proofs are rather tedious, in this section we only present outlines; readers can refer to Li and Zhu [18] for the details.

A.1. *Proofs of the theorems in Section 2.*

PROOF OF THEOREM 2.1. We first write out the formula for $J_n$. From definition (2.1), we have

$$\hat{\Sigma}(h) = \frac{1}{c(c-1)} \sum_{l=2}^{c} \sum_{j=1}^{l-1} (\boldsymbol{z}_{(h,l)} - \boldsymbol{z}_{(h,j)})^2.$$

For every $\boldsymbol{z}$, we have $\boldsymbol{z} = \mathbf{m}(y) + \boldsymbol{\varepsilon}$. Thus, for any pair $l$ and $j$,

$$
\begin{aligned}
&(\boldsymbol{z}_{(h,l)} - \boldsymbol{z}_{(h,j)})^2 \\
&\quad = (\mathbf{m}(y_{(h,l)}) - \mathbf{m}(y_{(h,j)}))^2 + (\mathbf{m}(y_{(h,l)}) - \mathbf{m}(y_{(h,j)}))(\boldsymbol{\varepsilon}_{(h,l)} - \boldsymbol{\varepsilon}_{(h,j)})^T \\
&\qquad + (\boldsymbol{\varepsilon}_{(h,l)} - \boldsymbol{\varepsilon}_{(h,j)})(\mathbf{m}(y_{(h,l)}) - \mathbf{m}(y_{(h,j)}))^T + (\boldsymbol{\varepsilon}_{(h,l)} - \boldsymbol{\varepsilon}_{(h,j)})^2 \\
&\quad =: \mathbf{S}_1(h,l,j) + \mathbf{S}_2(h,l,j) + \mathbf{S}_3(h,l,j) + \mathbf{S}_4(h,l,j).
\end{aligned}
$$

Further, $\Lambda_n$ can be written as

$$
\Lambda_n = \frac{\sum_{h=1}^{H}[\sum_{l=2}^{c}\sum_{j=1}^{l-1}(\mathbf{S}_1(h,l,j)+\mathbf{S}_2(h,l,j)+\mathbf{S}_3(h,l,j)+\mathbf{S}_4(h,l,j))]^2}{nc(c-1)^2}.
$$

For the sake of notational simplicity, we let

$$
\text{(A.1)} \quad C_n(i,k) = \frac{1}{nc(c-1)^2}\sum_{h=1}^{H}\sum_{l=2}^{c}\sum_{j=1}^{l-1}\sum_{v=2}^{c}\sum_{u=1}^{v-1}\mathbf{S}_i(h,l,j)\mathbf{S}_k(h,v,u).
$$

Then $\Lambda_n = \sum_{i=1}^{4}\sum_{k=1}^{4} C_n(i,k)$. Note that $A_n = C_n(4,4)$ and thus $J_n = \Lambda_n - C_n(4,4)$. To show that $n^\beta J_n = o_p(1)$, we only need to show that under the conditions of Theorem 2.1, for any pair $(i,k)$, except when $i = k = 4$, $n^\beta C_n(i,k)$ converges to 0 in probability as $n \to \infty$. Without loss of generality, we only consider the upper-left most element of $C_n(i,k)$, as the other elements can be handled similarly. Without confusion, we can still use the same notation for this element as the associated matrix $C_n(i,k)$. Therefore, in the following proof, $C_n(i,k)$ is real-valued.

For each $q$ such that $0 < q < \frac{1}{2}$, divide the outer summation over $h$ into three summations—from 1 to $[Hq]$, $[Hq]+1$ to $[H(1-q)]$ and $[H(1-q)]+1$ to $H$—to obtain

$$
C_n(i,k) = C_{1n}(i,k) + C_{2n}(i,k) + C_{3n}(i,k).
$$

For $C_{2n}(i,k)$, we have

$$
|C_{2n}(i,k)| \le \frac{1}{nc(c-1)^2}\sum_{h=[Hq]+1}^{[H(1-q)]}\sum_{l=1}^{c}\sum_{j=1}^{l-1}\sum_{v=2}^{c}\sum_{u=1}^{v-1}\|\mathbf{S}_i(h,l,j)\| \cdot \|\mathbf{S}_k(h,v,u)\|,
$$

where $\|\mathbf{S}\|$ denotes the maximum absolute value among elements in $\mathbf{S}$. For $\|\mathbf{S}_i(h,l,j)\| \cdot \|\mathbf{S}_k(h,v,u)\|$, we note that when $h \in [[Hq]+1, [H(1-q)]]$, there is a compact set $[-B(q), B(q)]$ such that in probability, both $y_{([nq]+1)}$ and $y_{([n(1-q)])}$ belong to that set. As $\mathbf{m}(y)$ is bounded on any compact set, there exists a $Q > 0$ such that in probability, $\|\mathbf{m}(y_{(h,j)})\| \le Q$. Let $\bar{\boldsymbol{\varepsilon}}_{(n)}$ and $\bar{\boldsymbol{\varepsilon}}_{(1)}$

denote the largest and the smallest of all $\varepsilon_{(i)}$'s, respectively. When $i$ and $k$ are fixed, we can determine $s$ such that

$$\sum_{l=2}^{c}\sum_{j=1}^{l-1}\sum_{v=2}^{c}\sum_{u=1}^{v-1}\|\mathbf{S}_i(h,l,j)\|\cdot\|\mathbf{S}_k(h,v,u)\|$$
$$\leq \frac{p^2c(c-1)\|\bar{\varepsilon}_{(n)}-\bar{\varepsilon}_{(1)}\|^{4-s}}{2}\sum_{l=2}^{c}\sum_{j=1}^{l-1}(2Q)^{s-1}\|\mathbf{m}(y_{(h,l)})-\mathbf{m}(y_{(h,j)})\|$$
$$+o_p(1).$$

As $i$ and $k$ cannot equal 4 simultaneously, we have $1\leq s\leq 4$ and hence,

$$\begin{aligned}C_{2n}(i,k)&\\ &\leq \frac{2^{s-2}\|\bar{\varepsilon}_{(n)}-\bar{\varepsilon}_{(1)}\|^{4-s}Q^{s-1}p^3c\sup_{\Pi_n(B(q))}\sum_{j=1}^{n-1}\|\mathbf{m}(y_{(j+1)})-\mathbf{m}(y_{(j)})\|}{n}\\ &\quad+o_p(1)\\ &=: C'_{2n}(s)+o_p(1).\end{aligned}$$

Using Lemma 1 of [14], we have $n^{-\frac{1}{8+\alpha}}\|\bar{\varepsilon}_{(1)}-\bar{\varepsilon}_{(1)}\|=o_p(1)$. Condition (2) of Theorem 2.1 implies that $\lim_{n\to\infty}n^{-r}\sup_{\Pi_n(B(q))}\sum_{i=1}^{n}\|\mathbf{m}(y_{(i+1)})-\mathbf{m}(y_{(i)})\|=0$. As $s\geq 1$, $C'_{2n}(s)=o_p(n^{r+\frac{3}{8+\alpha}+b-1})$ and therefore when $\beta+b+r+\frac{3}{8+\alpha}\leq 1$, $n^{\beta}C'_{2n}(s)\to 0$. We now consider $C_{1n}(i,k)$ and $C_{3n}(i,k)$. If $y$ is not bounded, we choose a sufficiently small $q$ so that $P(y_{([n(1-q)])}>B_0)\to 1$ as $n\to\infty$, where $B_0$ is given by condition (3) of Theorem 2.1. Using the nonexpansive property of $M(y)$, we can prove that

$$\begin{aligned}C_{3n}(i,k)&\leq \frac{p^3c\|\bar{\varepsilon}_{(n)}-\bar{\varepsilon}_{(1)}\|^{4-s}}{2n}\|M(y_{(n)})-M(y_{([n(1-q)])})\|^sI(y_{([n(1-q)])}>B_0)\\ &\quad+o_p(1)\\ &=: C'_{3n}(s)+o_p(1).\end{aligned}$$

By condition (3) and Lemma 1 of [14], it can be shown that when $\beta+b+\frac{4}{8+\alpha}\leq 1$, $n^{\beta}C'_{3n}(s)=o_p(1)$. The reasoning is similar for $C_{1n}(i,k)$, but we omit the details. The proof is thus complete. □

Proof of Theorem 2.2. The conditioning method is used to prove Theorem 2.2 and the other theorems. Denote $\mathcal{F}_n=\sigma\{y_1,\dots,y_n\}$. To compute $\mathbf{E}(A_n)$, we first compute the conditional expectation of $A_n$ given $y_i$'s as follows, where $A_n$ is defined in Section 2.1:

$$\mathbf{E}(A_n|\mathcal{F}_n)$$

$$
\begin{aligned}
&= \sum_{h=1}^{H}\sum_{l=1}^{c}\frac{\mathbf{E}((\boldsymbol{\varepsilon}_{(h,l)}\boldsymbol{\varepsilon}_{(h,l)}^T)^2|\mathcal{F}_n)}{nc} \\
&\quad + \sum_{h=1}^{H}\sum_{l=1}^{c}\sum_{v=1(v\neq l)}^{c}\frac{1}{nc}\left(1+\frac{1}{(c-1)^2}\right)\mathbf{E}(\boldsymbol{\varepsilon}_{(h,l)}\boldsymbol{\varepsilon}_{(h,l)}^T|\mathcal{F}_n)\mathbf{E}(\boldsymbol{\varepsilon}_{(h,v)}\boldsymbol{\varepsilon}_{(h,v)}^T|\mathcal{F}_n) \qquad \text{(A.2)} \\
&\quad + \sum_{h=1}^{H}\sum_{l=1}^{c}\sum_{v=1(v\neq l)}^{c}\frac{1}{nc(c-1)^2}\mathbf{E}((\boldsymbol{\varepsilon}_{(h,l)}\boldsymbol{\varepsilon}_{(h,v)}^T)^2|\mathcal{F}_n) \\
&=: \mathbf{E}(A_{1n}|\mathcal{F}_n)+\mathbf{E}(A_{2n}|\mathcal{F}_n)+\mathbf{E}(A_{3n}|\mathcal{F}_n).
\end{aligned}
$$

As the $\boldsymbol{\varepsilon}_{(i)}$'s are conditionally independent when the $y_i$ are given, $\mathbf{E}(A_{1n}|\mathcal{F}_n)$ is equal to $\frac{1}{nc}\sum_{j=1}^n \mathbf{E}((\boldsymbol{\varepsilon}_j\boldsymbol{\varepsilon}_j^T)^2|y_j)$. This is a sum of i.i.d. random variables and therefore $\mathbf{E}(A_{1n}) = \frac{1}{c}\mathbf{E}[(\boldsymbol{\varepsilon}\boldsymbol{\varepsilon}^T)^2]$. For $\mathbf{E}(A_{2n}|\mathcal{F}_n)$, the conditional independence property and the definition $\mathbf{m}_1(y)=\mathbf{E}(\boldsymbol{\varepsilon}\boldsymbol{\varepsilon}^T|y)$ together yield that

$$
\begin{aligned}
&\mathbf{E}(A_{2n}|\mathcal{F}_n) \\
&\quad = \frac{(c-1)((c-1)^2+1)}{nc(c-1)^2}\sum_{h=1}^{H}\sum_{l=1}^{c}\mathbf{m}_1(y_{(h,l)})\mathbf{m}_1(y_{(h,l)})^T \\
&\qquad + \frac{(c-1)^2+1}{nc(c-1)^2}\sum_{h=1}^{H}\sum_{l=1}^{c}\sum_{v=1(v\neq l)}^{c}\mathbf{m}_1(y_{(h,l)})(\mathbf{m}_1(y_{(h,v)})-\mathbf{m}_1(y_{(h,l)}))^T \\
&\quad =: \mathbf{E}(A_{21n}|\mathcal{F}_n)+\mathbf{E}(A_{22n}|\mathcal{F}_n).
\end{aligned}
$$

As $\mathbf{E}(A_{21n}|\mathcal{F}_n) = \frac{1}{n}(1-\frac{(c-2)}{c(c-1)})\sum_{j=1}^n \mathbf{m}_1(y_j)^2$, we have that $\mathbf{E}(A_{21n}) = (1-\frac{(c-2)}{c(c-1)})\Lambda$.

For $\mathbf{E}(A_{22n}|\mathcal{F}_n)$, the conclusion is

$$
\mathbf{E}(A_{22n}|\mathcal{F}_n) = o_p(cn^{-1+\max\{r_1,\frac{2}{4+\alpha/2}\}}). \qquad \text{(A.3)}
$$

The lines of the proof essentially follow those of the proof of Theorem 2.1. For each $q_1$ such that $0<q_1<\frac{1}{2}$, we divide the outer summation over $h$ into three summations: from 1 to $[Hq_1]$, $[Hq_1]+1$ to $[H(1-q_1)]$ and $[H(1-q_1)]+1$ to $H$. Hence, $\mathbf{E}(A_{22n}|\mathcal{F}_n) = D_{1n}+D_{2n}+D_{3n}$. Note that when $h\in[[Hq_1]+1,[H(1-q_1)]]$, there exists a constant $Q_1$ such that $\|\mathbf{m}_1(y_{(h,l)})\|\le Q_1$ for all $1\le l\le c$. Thus, as $\mathbf{m}_1(y)$ has total variation of order $r_1$,

$$
\begin{aligned}
D_{2n} &\le \frac{Q_1((c-1)^2+1)p^3\sup_{\Pi_n(B(q_1))}\sum_{i=1}^n\|\mathbf{m}_1(y_{(i+1)})-\mathbf{m}_1(y_{(i)})\|}{n(c-1)} + o_p(1) \\
&= o(cn^{-1+r_1}).
\end{aligned}
$$

If $y$ is not bounded, then we choose a sufficiently small $q_1$ so that $P(y_{([n(1-q_1)])} > B_0') \to 1$ as $n \to \infty$, where $B_0'$ is given by condition (3) of Theorem 2.2. Similarly, $D_{3n} = o_p(cn^{-1+\frac{2}{4+\alpha/2}})$. The proof is similar to that for $D_{1n}$ and (A.3) then holds. By condition (5) and Lemma 4.11 of [15], we have

$$\mathbf{E}(A_{22n}) = o(cn^{-1+\max\{r_1, \frac{2}{4+\alpha/2}, \rho_1\}}). \tag{A.4}$$

The proof of $\mathbf{E}(A_{3n}|\mathcal{F}_n)$ of (A.2) is very similar to the one just given and we can thus obtain $\mathbf{E}(A_{3n}) = o(c^{-1}n^{-1+\max\{r_1, \frac{2}{4+2/\alpha}, \rho_1\}})$. Hence, (2.3) is proved.

We now turn to the proof of the second conclusion, (2.4), that $n^\beta(A_n - \Lambda) = o_p(1)$. Without loss of generality, consider the upper-rightmost element of $n^\beta(A_n - \Lambda)$. Without confusion, we can still use the notation $n^\beta(A_n - \Lambda)$ to represent this element. Note that $n^\beta\{A_n - \Lambda\} = n^\beta\{A_n - \mathbf{E}(A_n|\mathcal{F}_n) + \mathbf{E}(A_n|\mathcal{F}_n) - \Lambda\}$. From the proof of (2.3), we can obtain that when $\beta < b$ and $\beta \le 1 - b - \max\{r_1, \frac{2}{4+\alpha/2}\}$,

$$n^\beta\{\mathbf{E}(A_n|\mathcal{F}_n) - \Lambda\} = o_p(1). \tag{A.5}$$

Therefore, it remains to show that $n^\beta\{A_n - \mathbf{E}(A_n|\mathcal{F}_n)\} = o_p(1)$ and it suffices to demonstrate the convergence of its second moment. That is, as $n \to \infty$,

$$n^{2\beta}\mathbf{E}[(\{(A_n - \mathbf{E}(A_n|\mathcal{F}_n))\})^2] \to 0. \tag{A.6}$$

Invoking (A.2), the definition of $A_n$ given in Section 2.1, and rearranging the terms, we see that

$$\begin{aligned}
&(A_n - \mathbf{E}(A_n|\mathcal{F}_n)) \\
&\quad = \frac{1}{n}\sum_{h=1}^{H}\Bigg\{\Bigg[\frac{1}{c}\sum_{l=1}^{c}\sum_{v=1(v\neq l)}^{c} \varepsilon_{(h,l)}^2\varepsilon_{(h,v)}^2 \\
&\qquad\qquad - \frac{1}{c}\sum_{l=1}^{c}\sum_{v=1(v\neq l)}^{c}(\mathbf{E}(\varepsilon_{(h,l)}^2|y_{(h,l)}))(\mathbf{E}(\varepsilon_{(h,v)}^2|y_{(h,v)}))\Bigg] \\
&\qquad\qquad + \Bigg[\frac{1}{c}\sum_{l=1}^{c}((\varepsilon_{(h,l)}\varepsilon_{(h,l)}^T)^2 - \mathbf{E}((\varepsilon_{(h,l)}\varepsilon_{(h,l)}^T)^2|y_{(h,l)}))\Bigg] \\
&\qquad\qquad + \Bigg[\frac{1}{c(c-1)^2}\sum_{l=1}^{c}\sum_{j=1(j\neq l)}^{c}\sum_{v=1}^{c}\sum_{u=1(u\neq v)}^{c}\varepsilon_{(h,l)}\varepsilon_{(h,j)}^T\varepsilon_{(h,v)}\varepsilon_{(h,u)}^T \\
&\qquad\qquad - \frac{1}{c(c-1)^2}\sum_{l=1}^{c}\sum_{v=1(v\neq l)}^{c}(\mathbf{E}(\varepsilon_{(h,l)}^2|y_{(h,l)}))(\mathbf{E}(\varepsilon_{(h,v)}^2|y_{(h,v)}))\Bigg]
\end{aligned}$$

$$- \left[ \frac{1}{c(c-1)} \left( \sum_{l=1}^{c} \sum_{v=1}^{c} \sum_{u=1(u \neq v)}^{c} \boldsymbol{\varepsilon}_{(h,l)}^2 \boldsymbol{\varepsilon}_{(h,v)} \boldsymbol{\varepsilon}_{(h,u)}^T + \sum_{l=1}^{c} \sum_{j=1(j \neq l)}^{c} \sum_{v=1}^{c} \boldsymbol{\varepsilon}_{(h,l)} \boldsymbol{\varepsilon}_{(h,j)}^T \boldsymbol{\varepsilon}_{(h,v)}^2 \right) \right] \Bigg\}$$

$$=: \frac{1}{n} \sum_{h=1}^{H} \{V_0(h) + V_1(h) + V_2(h) + V_3(h)\}.$$

We again use the conditioning method to show that $\frac{n^{2\beta}}{n^2} \sum_{h=1}^{H} \mathbf{E} V_i^2(h) = o(1)$ for $i = 0$, 1, 2 and 3 and then use the inequality $2|V_i(h)V_j(h)| \le V_i^2(h) + V_j^2(h)$ to obtain that the intersection terms converge to zero from the convergence of $\mathbf{E}(V_i^2(h))$. The proof of Theorem 2.2 can then be completed. We now proceed to the first step as follows.

To simplify the notation, we write, for any integer $l > 1$, $\mathbf{E}^l(\boldsymbol{\varepsilon}^s|y) = \mathbf{E}^{l-1}(\boldsymbol{\varepsilon}^s|y)\mathbf{E}(\boldsymbol{\varepsilon}^s|y)$, where $1 \le s \le 6$. By means of elementary calculation, we obtain the result

$$\frac{n^{2\beta}}{n^2} \sum_{h=1}^{H} \mathbf{E}(V_1^2(h)) = O\left( \frac{n^{2\beta}}{nc^2} \mathbf{E} \boldsymbol{\varepsilon}^8 - \frac{n^{2\beta}}{nc^2} \mathbf{E}(\mathbf{E}^2(\boldsymbol{\varepsilon}^4|y)) \right) = o(1).$$

$\frac{n^{2\beta}}{n^2} \sum_{h=1}^{H} \mathbf{E}(V_2^2(h))$ can be bounded by

$$\left( \frac{56 n^{2\beta}}{nc^3} \mathbf{E}(\mathbf{E}^4(\boldsymbol{\varepsilon}^2|y)) + \frac{64 n^{2\beta}}{nc^4} \mathbf{E}(\mathbf{E}^3(\boldsymbol{\varepsilon}^3|y)) + \frac{16 n^{2\beta}}{nc^4} \mathbf{E}(\mathbf{E}^3(\boldsymbol{\varepsilon}^4|y)) + \frac{64 n^{2\beta}}{nc^4} \mathbf{E}(\mathbf{E}^3(\boldsymbol{\varepsilon}^2|y)) + \frac{8 n^{2\beta}}{nc^5} \mathbf{E}\mathbf{E}^2(\boldsymbol{\varepsilon}^4|y) \right).$$

Since $\mathbf{E}(\boldsymbol{\varepsilon}^{12}) < \infty$, it is $o_p(1)$. Similarly, we have $\frac{n^{2\beta}}{n^2} \sum_{h=1}^{H} \mathbf{E}(V_3^2(h)) = o_p(1)$.

Using the conditioning method, we can also prove that the sum that relates to $\mathbf{E}(V_0^2(h))$ converges to zero. First, we have

$$\begin{aligned}
\mathbf{E}(V_0^2(h)|\mathcal{F}) = & \left[ \frac{2}{c^2} \sum_{l=1}^{c} \sum_{j=1(l \neq j)}^{c} \mathbf{E}(\boldsymbol{\varepsilon}_{(h,l)}^4|\mathcal{F}) \mathbf{E}(\boldsymbol{\varepsilon}_{(h,j)}^4|y_{(h,j)}) \right] \\
& - \left[ \frac{2}{c^2} \sum_{l=1}^{c} \sum_{j=1(l \neq j)}^{c} \mathbf{E}^2(\boldsymbol{\varepsilon}_{(h,l)}^2|y_{(h,l)}) \mathbf{E}^2(\boldsymbol{\varepsilon}_{(h,j)}^2|y_{(h,j)}) \right] \\
& + \left[ \frac{4c^2}{c^2} \sum_{l=1}^{c} (\mathbf{E}(\boldsymbol{\varepsilon}_{(h,l)}^4|y_{(h,l)}) \mathbf{E}^2(\boldsymbol{\varepsilon}_{(h,l)}^2|y_{(h,l)}) - \mathbf{E}^4(\boldsymbol{\varepsilon}_{(h,l)}^2|y_{(h,l)})) \right] \\
& - \left[ \frac{4}{c^2} \sum\sum\sum_{1 \le l \neq j \neq v \le c} u_{h,l,j,v}^1 \right] - \left[ \frac{4}{c^2} \sum\sum\sum_{1 \le l \neq j \neq u \le c} u_{h,l,j,v}^2 \right]
\end{aligned}$$

$$+\left[\frac{4}{c^2}\sum\sum_{1\le l\ne j\ne u\le c}\sum u^3_{h,l,j,v}\right]+\left[\frac{4}{c^2}\sum\sum_{1\le l\ne j\ne u\le c}\sum u^4_{h,l,j,v}\right]$$
$$=: V_{00}(h)-V_{01}(h)+V_{02}(h)-V_{03}(h)-V_{04}(h)$$
$$+V_{05}(h)+V_{06}(h),$$

where

$$u^1_{h,l,j,v}=\mathbf{m}_2(y_{(h,l)})(\mathbf{m}_1(y_{(h,l)})-\mathbf{m}_1(y_{(h,v)}))\mathbf{m}_1(y_{(h,l)}),$$
$$u^2_{h,l,j,v}=\mathbf{m}_2(y_{(h,l)})\mathbf{m}_1(y_{(h,v)})(\mathbf{m}_1(y_{(h,l)})-\mathbf{m}_1(y_{(h,j)})),$$
$$u^3_{h,l,j,v}=\mathbf{m}^2_1(y_{(h,l)})(\mathbf{m}_1(y_{(h,l)})-\mathbf{m}_1(y_{(h,v)}))\mathbf{m}_1(y_{(h,l)}),$$
$$u^4_{h,l,j,v}=\mathbf{m}^2_1(y_{(h,l)})\mathbf{m}_1(y_{(h,v)})(\mathbf{m}_1(y_{(h,l)})-\mathbf{m}_1(y_{(h,j)})).$$

We now prove that when $c\sim n^b$ and $2\beta+\max\{2r_1,\frac{1}{2+\alpha/4}+\frac{2}{4+\alpha/2},\rho_2\}+b\le 2$, all of the terms $\frac{n^{2\beta}}{n^2}\sum_{h=1}^H\mathbf{E}(V_{0i}(h))$ tend to 0. Using the conditioning method and the inequality

$$\mathbf{E}(\varepsilon^4_{(h,l)}|y_{(h,l)})\mathbf{E}(\varepsilon^4_{(h,j)}|y_{(h,j)})\le\frac{1}{2}(\mathbf{E}^2(\varepsilon^4_{(h,l)}|y_{(h,l)})+\mathbf{E}^2(\varepsilon^4_{(h,j)}|y_{(h,j)})),$$

we have

$$\frac{n^{2\beta}}{n^2}\sum_{h=1}^H\mathbf{E}V_{00}(h)=O\biggl(\frac{2n^{2\beta}}{nc}\mathbf{E}(\mathbf{E}^2(\varepsilon|y))\biggr)=o(1).$$

Similar arguments can be used to obtain $\frac{n^{2\beta}}{n^2}\sum_{h=1}^H\mathbf{E}(V_{01}(h))=o(1)$.

As $V_{02}(h)$ is a sum of i.i.d. random variables, invoking the conditions of Theorem 2.2, the fact that $\beta<0.5$ and the law of large numbers, we can show that $\frac{n^{2\beta}}{n^2}\sum_{h=1}^H V_{02}(h)=o(1)$.

The proof of the sum of $V_{03}(h)$ is similar to that of $\mathbf{E}(A_{22n}|\mathcal{F}_n)$. We choose $0<q_2<1$ and divide the summation of $h$ into three parts: $[1,[Hq_2]]$, $[[Hq_2]+1,[H(1-q_2)]]$ and $[[H(1-q_2)]+1,H]$. The sums of the conditional expectation of $\mathbf{E}(V_{03}(h)|\mathcal{F}_n)$ over $h$ in these three intervals are analyzed and $\frac{n^{2\beta}}{n^2}\sum_{h=[Hq_2]+1}^{[H(1-q_2)]}\mathbf{E}(V_{03}(h))$ can be proved to be asymptotically zero. The proof is very similar to that of (A.3) and thus we omit the details in this paper. The proof of (2.4) is thus complete.

This completes proof of Theorem 2.2. □

PROOF OF THEOREM 2.3. The proof is similar to that of Theorem 3.1 below, and thus we omit the details. □

A.2. *Proofs of the theorems in Section 3.*

PROOF OF THEOREM 3.1. Our goal is to determine the asymptotic behavior of $\frac{1}{H}\sum_{h=1}^{H}(I_p-\hat{\Sigma}(h))^2$, where $\hat{\Sigma}(h)$ is defined in (3.1) and $S_h=(y_{(c(h-1))},y_{(ch)}]$. It suffices to show that for any $p(p+1)/2$ vector $\mathbf{a}$, $\mathbf{a}^T\operatorname{vech}\{\frac{1}{H}\sum_{h=1}^{H}(I_p-\hat{\Sigma}(h))^2\}$ is asymptotically univariate normal. Again, for the sake of notational simplicity, we consider the univariate case. Clearly, $\hat{q}_h=y_{(ch)}$, $h=1,\ldots,H$, are the empirical quantiles that converge to the population quantiles $q_h$ in probability, where $P(Y\le q_h)=h/H$. If we can verify the asymptotic normality of $\hat{\Sigma}(h)-\Sigma(h)$ for $h=1,\ldots,H$, then the asymptotic normality of $\Lambda_n$ can be obtained through the decomposition

$$
\begin{aligned}
&\sqrt{n}\left(\frac{1}{H}\sum_{h=1}^{H}(I_p-\hat{\Sigma}(h))^2-\frac{1}{H}\sum_{h=1}^{H}((I_p-\Sigma(h))^2\right)\\
&\quad=\frac{-\sqrt{n}}{H}\sum_{h\ =1}^{H}(\hat{\Sigma}(h)-\Sigma(h))(2I_p-\hat{\Sigma}(h)-\Sigma(h))\\
&\quad=\frac{-2\sqrt{n}}{H}\sum_{h\ =1}^{H}(\hat{\Sigma}(h)-\Sigma(h))(I_p-\Sigma(h))+o_p(1).
\end{aligned}\tag{A.7}
$$

We now study $\hat{\Sigma}(h)$. From (3.1),

$$
\begin{aligned}
\hat{\Sigma}(h)&=\frac{1}{n\hat{p}_h}\sum_{j=1}^{n}\boldsymbol{z}_j^2I(y_j\in\hat{S}_h)-\left(\frac{1}{n\hat{p}_h}\sum_{j=1}^{n}\boldsymbol{z}_jI(y_j\in\hat{S}_h)\right)^2\\
&=\hat{\Sigma}_1(h)-(\hat{E}_1(h))^2.
\end{aligned}\tag{A.8}
$$

Next, we calculate $\sqrt{n}(\hat{\Sigma}_1(h)-\Sigma_1(h))$. Note that $\hat{p}_h=p_h=1/H$ and thus

$$
\begin{aligned}
\sqrt{n}(\hat{\Sigma}_1(h)-\Sigma_1(h))&=\frac{1}{\sqrt{n}\hat{p}_h}\sum_{j=1}^{n}(\boldsymbol{z}_j^2I(y_j\in S_h)-\mathbf{E}(\boldsymbol{z}_j^2I(y_j\in S_h)))\\
&\quad+\frac{1}{\sqrt{n}\hat{p}_h}\sum_{j=1}^{n}\boldsymbol{z}_j^2(I(y_j\in\hat{S}_h)-I(y_j\in S_h))\\
&=:\hat{\Sigma}_{11}(h)+\hat{\Sigma}_{12}(h).
\end{aligned}\tag{A.9}
$$

Clearly, $\hat{\Sigma}_{11}(h)$ is asymptotically normal because it is a sum of i.i.d. random variables.

For $\hat{\Sigma}_{12}(h)$, we first introduce the notation $F(Y,\boldsymbol{z},a,b)=\boldsymbol{z}^2(I(Y\in(a,b])-I(Y\in S_h))$ for any pair $(a,b)$. Note that $\hat{q}_h-q_h=O_p(1/\sqrt{n})$. Invoking Theorem 1 of Zhu and Ng [27] or the argument used in Stute and Zhu [22] and

Stute, Thies and Zhu [21], we can show that

$$\left|\frac{1}{\sqrt{n p_h}}\sum_{j=1}^{n}(F(y_j,\boldsymbol{z}_j,\hat{q}_{h-1},\hat{q}_h)-\mathbf{E}(F(Y,\boldsymbol{z},\hat{q}_{h-1},\hat{q}_h)))\right|=o_p(1).$$

Together with (A.10), the continuity of $\mathbf{E}(F(Y,\boldsymbol{z},q_{h-1},q_h))$ at $q_{h-1}$ and $q_h$, the $\sqrt{n}$ consistency of $q_h$ and Taylor expansion give

$$\begin{aligned}\hat{\Sigma}_{12}(h)&=H\sqrt{n}\mathbf{E}(F(Y,\boldsymbol{z},\hat{q}_{h-1},\hat{q}_h))+o_p(1)\\&=H\sqrt{n}(\hat{q}_{h-1}-q_{h-1},\hat{q}_h-q_h)\tilde{F}'(q_{h-1},q_h)+o_p(1)\\&=\frac{H}{\sqrt{n}}\sum_{j=1}^{n}\left(\frac{-I(y_j\le q_{h-1})+\frac{h-1}{H}}{f(q_{h-1})},\frac{-I(y_j\le q_h)+\frac{h}{H}}{f(q_h)}\right)\tilde{F}'(q_{h-1},q_h)\\&\quad+o_p(1),\end{aligned}\tag{A.10}$$

where $\tilde{F}'$ is the derivative of $\mathbf{E}(F(Y,\boldsymbol{z},a,b))$ with respect to $(a,b)$. The asymptotic normality can be shown to hold by using well-known results on the empirical quantiles $\hat{q}_h$ (see [20]).

For $(\hat{E}_1(h))^2$ from (A.8), the foregoing argument can be applied to obtain $\sqrt{n}(\hat{E}_1(h))^2$, giving

$$\begin{aligned}&\sqrt{n}((\hat{E}_1(h))^2-(E_1(h))^2)\\&\quad=2\sqrt{n}(\hat{E}_1(h)-E_1(h))E_1(h)+o_p(1)\\&\quad=\frac{2H}{\sqrt{n}}\sum_{j=1}^{n}(\boldsymbol{z}_jI(y_j\in S_h)-\mathbf{E}(\boldsymbol{z}I(Y\in S_h)))E_1(h)\\&\qquad+\frac{2H}{\sqrt{n}}\sum_{j=1}^{n}\left(\frac{-I(y_j\le q_{h-1})+\frac{h-1}{H}}{f(q_{h-1})},\frac{-I(y_j\le q_h)+\frac{h}{H}}{f(q_h)}\right)\\&\qquad\times\tilde{G}'(q_{h-1},q_h)E_1(h)+o_p(1),\end{aligned}\tag{A.11}$$

where $\tilde{G}'(a,b)$ is the derivative of $\mathbf{E}(G(Y,\boldsymbol{z},a,b)):=\mathbf{E}(\boldsymbol{z}(I(Y\in(a,b])-I(Y\in S_h)))$ with respect to $(a,b)$. Together with (A.8)–(A.12), we have

$$\begin{aligned}&\sqrt{n}\left(\frac{1}{H}\sum_{h=1}^{H}(I_p-\hat{\Sigma}(h))^2-\frac{1}{H}\sum_{h=1}^{H}(I_p-\Sigma(h))^2\right)\\&\quad=\frac{1}{\sqrt{n}}\sum_{j=1}^{n}\left\{-2\sum_{h=1}^{H}((\boldsymbol{z}_j^2-2\boldsymbol{z}_jE_1(h))I(y_j\in S_h)\right.\end{aligned}$$

$$
\begin{aligned}
&\qquad - \mathbf{E}((\boldsymbol{z}^2 - 2\boldsymbol{z}E_1(h))I(Y \in S_h))) \\
&\qquad - 2\sum_{h=1}^{H}\left(\frac{-I(y_j \le q_{h-1}) + \frac{h-1}{H}}{f(q_{h-1})}, -\frac{I(y_j \le q_h) + \frac{h}{H}}{f(q_h)}\right) \\
&\qquad \times (\tilde{F}'(q_{h-1}, q_h) - 2\tilde{G}'(q_{h-1}, q_h)E_1(h))\Bigg\} \times (I_p - \Sigma(h)) \\
&+ o_p(1) \\
&:= \frac{1}{\sqrt{n}}\sum_{j=1}^{n} L(y_j, z_j) + o_p(1) \Rightarrow N(0, \Delta'),
\end{aligned}
$$

where $\Delta' = \mathrm{Cov}(L(Y, \boldsymbol{z}))$. $\square$

Proof of Theorem 3.2. We only present the proof for the univariate case. As $c \to \infty$, it is equivalent to showing that when $c$ satisfies the required conditions,

$$
\frac{\sqrt{n}}{c}\left(\frac{1}{H}\sum_{h=1}^{H}\frac{1}{c}\sum_{j=1}^{c}(\boldsymbol{z}_{(h,j)} - \bar{\boldsymbol{z}}_{(h)})^4 - \mathbf{E}(\boldsymbol{\varepsilon}^4)\right) = o_p(1). \tag{A.12}
$$

Some elementary calculation yields

$$
\begin{aligned}
&\frac{1}{H}\sum_{h=1}^{H}\frac{1}{c}\sum_{j=1}^{c}(\boldsymbol{z}_{(h,j)} - \bar{\boldsymbol{z}}_{(h)})^4 \\
&\quad = \frac{1}{H}\sum_{h=1}^{H}\frac{1}{c}\sum_{j=1}^{c}\boldsymbol{\varepsilon}_{(h,j)}^4 \\
&\qquad + \frac{1}{H}\sum_{h=1}^{H}\frac{1}{c}\sum_{j=1}^{c}\left(\frac{-4\boldsymbol{\varepsilon}_{(h,j)}^3}{c}(A_{(h)} + B_{(h,j)}) + \frac{6\boldsymbol{\varepsilon}_{(h,j)}^2}{c^2}(A_{(h)} + B_{(h,j)})^2\right. \qquad \text{(A.13)} \\
&\qquad\qquad \left. - \frac{4\boldsymbol{\varepsilon}_{(h,j)}}{c^3}(A_{(h)} + B_{(h,j)})^3 + \frac{1}{c^4}(A_{(h)} + B_{(h,j)})^4\right) \\
&\quad =: R_{n1} + R_{n2},
\end{aligned}
$$

where $A_{(h)} = \sum_{v=1}^{c}\boldsymbol{\varepsilon}_{(h,v)}$ and $B_{(h,j)} = \sum_{v=1}^{c}(\mathbf{m}(y_{(h,v)}) - \mathbf{m}(y_{(h,j)}))$. Rearranging the summands in $R_{n1}$, we can easily show that $\sqrt{n}[R_{n1} - \mathbf{E}(\boldsymbol{\varepsilon}^4)] = \frac{1}{\sqrt{n}}\sum_{j=1}^{n}(\boldsymbol{\varepsilon}_j^4 - \mathbf{E}(\boldsymbol{\varepsilon}^4))$ follows the distribution $N(0, \mathrm{var}(\boldsymbol{\varepsilon}^4))$ and thus $\frac{\sqrt{n}}{c}[R_{n1} - \mathbf{E}(\boldsymbol{\varepsilon}^4)] = o_p(1)$. Hence, to prove (A.12), we only need to show that

$$
\frac{\sqrt{n}}{c}R_{n2} = o_p(1). \tag{A.14}
$$

We find that the terms in $\frac{\sqrt{n}}{c}R_{n2}$ have the following two common formats. For $1 \le s_1 \le 4$,

$$K(s_1) := \frac{\sqrt{n}}{c}\frac{1}{H}\sum_{h=1}^{H}\frac{1}{c}\sum_{j=1}^{c}\boldsymbol{\varepsilon}_{(h,j)}^{4-s_1}\frac{1}{c^{s_1}}A_{(h)}^{s_1}, \tag{A.15}$$

and for $1 \le s' \le 4$ and $0 \le s \le 4 - s'$,

$$W(s,s') := \frac{\sqrt{n}}{c}\frac{1}{H}\sum_{h=1}^{H}\frac{1}{c}\sum_{j=1}^{c}\boldsymbol{\varepsilon}_{(h,j)}^{s}\frac{1}{c^{4-s}}A_{(h)}^{4-s-s'}B_{(h,j)}^{s'}. \tag{A.16}$$

Therefore, our task is to prove that they are all $o_p(1)$. For $K(s_1)$'s, we need only show that their second moments asymptotically converge to 0, the main idea of which is to use the conditioning method to compute their conditional expectations given $y_i$'s and to use a sum of i.i.d. random variables to approximate the $K(s_1)$'s. The arguments are very similar to those in the proof of Theorem 2.1 and the details can be found in [18].

For $W(s,s')$ of (A.16), we note that if we let $d = \max_{1\le i\le n}(|\boldsymbol{\varepsilon}_i|)$, then $|\frac{A_{(h)}}{c}| \le d$ and thus

$$W(s,s') \le \frac{\sqrt{n}d^{4-s'}}{c^{2+s'}}\frac{1}{H}\sum_{h=1}^{H}\sum_{j=1}^{c}B_{(h,j)}^{s'}.$$

For each $q$ such that $0 < q < \frac{1}{2}$, we divide the outer summation over $h$ into three summations—from 1 to $[Hq]$, $[Hq]+1$ to $[H(1-q)]$ and $[H(1-q)]+1$ to $H$—which allows us to write $W(s,s') = W_1(s,s') + W_2(s,s') + W_3(s,s')$. We then use the argument that was used to prove Theorem 2.1 to show that $W(s,s') = o_p(1)$. (A.14) is thus proved and the proof of Theorem 3.2 is complete. $\square$

Proof of Corollary 3.1. We want to show that for any $p(p+1)/2$ vector $\mathbf{a}$, $\mathbf{a}^T\mathrm{vech}\{\mathrm{CSAVE}_n - \Lambda\}$ is asymptotically univariate normal with zero mean and finite variance. Denote

$$Z_{nh} = \mathbf{a}^T\mathrm{vech}\Bigg\{\frac{(c-1)}{(c-1)^2+1}\sum_{l=1}^{c}\sum_{v=1}^{c}(\boldsymbol{\varepsilon}_{(h,l)}^2\boldsymbol{\varepsilon}_{(h,v)}^2) - c\Lambda - \frac{1}{c}\sum_{j=1}^{c}(\boldsymbol{\varepsilon}_{(h,j)} - \bar{\boldsymbol{\varepsilon}}_{(h)})^4$$
$$- \frac{2}{c-1}\sum_{l=2}^{c}\sum_{j=1}^{l-1}((\boldsymbol{\varepsilon}_{(h,l)} - \boldsymbol{\varepsilon}_{(h,j)})^2 - 2\mathbf{E}(\Sigma_{\boldsymbol{z}|y}))\Bigg\}.$$

To prove the asymptotic normality, we will check the four conditions with the conditional central limit theorem (CCLT) that was provided by Hsing and Carroll [14], Theorem A.4. From Theorem 3.2, $\sqrt{n}\mathbf{a}^T\mathrm{vech}\{\mathrm{CSAVE}_n - \mathbf{E}(I_p - \Sigma_{\boldsymbol{z}|y})^2\}$ is asymptotically equivalent to $\frac{1}{\sqrt{n}}\sum_{h=1}^{H}Z_{nh}$. As $Z_{n1},\dots,Z_{nH}$

are conditionally independent given $\mathcal{F}_n$, condition (1) of the CCLT is satisfied.

To check conditions (2)–(4) of the CCLT, the calculation is very similar to that in the proofs of Theorem 2.2 and Theorem 3.2. For the conditional expectation of $Z_{nh}$, we have

$$\begin{aligned}
&\frac{1}{\sqrt{n}}\sum_{h=1}^{H}\mathbf{E}(Z_{nh}|\mathcal{F}_n)\\
&= \frac{1}{\sqrt{n}}\sum_{j=1}^{n}\mathbf{a}^T\text{vech}\{\mathbf{m}_1^2(y_{(j)}) - \Lambda - 2(\mathbf{m}_1(y_{(j)}) - \mathbf{E}(\Sigma_{\boldsymbol{z}|y}))\} + o_p(1) \qquad \text{(A.17)}\\
&\to_d N(0, \mathbf{a}^T\Delta_1\mathbf{a}),
\end{aligned}$$

where $\Delta_1 = \text{var}(\text{vech}\{\mathbf{m}_1^2(y_{(j)}) - \Lambda - 2(\mathbf{m}_1(y_{(j)}) - \mathbf{E}(\Sigma_{\boldsymbol{z}|y}))\})$, and hence condition (4) of the CCLT is satisfied. For condition (2), we only need to note that, together with conditional independence,

$$\begin{aligned}
&\frac{1}{n}\sum_{h=1}^{H}\mathbf{E}\{(Z_{nh} - \mathbf{E}(Z_{nh}|\mathcal{F}_n))^2|\mathcal{F}_n\}\\
&\quad = \frac{1}{n}\sum_{j=1}^{n}\mathbf{a}^T\text{vech}\{(\mathbf{m}_2(y_{(j)}) - \mathbf{m}_1^2(y_{(j)}))\mathbf{m}_1^2(y_{(j)})\}\mathbf{a}\\
&\qquad + \frac{4}{n}\sum_{j=1}^{n}\mathbf{a}^T\text{vech}\{\mathbf{m}_2(y_{(j)}) - \mathbf{m}_1^2(y_{(j)})\}\mathbf{a}\\
&\qquad - \frac{4}{n}\sum_{j=1}^{n}\mathbf{a}^T\text{vech}\{(\mathbf{m}_2(y_{(j)}) - \mathbf{m}_1^2(y_{(j)}))\mathbf{m}_1(y_{(j)})\}\mathbf{a} + o_p(1) \qquad \text{(A.18)}\\
&\quad = \mathbf{a}^T\text{vech}\{\mathbf{E}[(\mathbf{m}_2(y) - \mathbf{m}_1^2(y))\mathbf{m}_1^2(y) + 4(\mathbf{m}_2(y) - \mathbf{m}_1^2(y))\\
&\qquad\qquad - 4(\mathbf{m}_2(y) - \mathbf{m}_1^2(y))\mathbf{m}_1(y)]\}\mathbf{a} + o_p(1)\\
&\quad =: \mathbf{a}^T\Delta_2\mathbf{a} + o_p(1).
\end{aligned}$$

Condition (3) of the CCLT can be checked using a similar argument. The main idea is as follows. Invoking the conditional independence of the $Z_{nh}$'s and the existence of the 12th moment, we can use a method similar to that which was used to prove Liapounoff's central limit theorem (see, e.g., Pollard [19]) to verify condition (3) of the CCLT. Hence, the CCLT implies that $\frac{1}{\sqrt{n}}\sum_{h=1}^{H} Z_{nh}$ is asymptotically normal with zero mean and variance $\mathbf{a}^T(\Delta_1 + \Delta_2)\mathbf{a}$.

When the $\hat{\boldsymbol{z}}_i$'s are used to construct the statistic, as with the proofs of the other theorems, the asymptotic normality holds with limiting variance $\mathbf{a}^T(\Delta_1 + \Delta_2 + \mathbf{E}_1)\mathbf{a}$, where $\mathbf{E}_1$ is the random matrix defined in (2.11). The proof is thus complete. $\square$

**Acknowledgment.** The first version of this paper was written when the two authors were at the University of Hong Kong.

Department of Statistics
Cornell University
Ithaca, New York 14853
USA
E-mail: yl377@cornell.edu

Department of Mathematics
Hong Kong Baptist University
Kowloon Tong
Hong Kong
E-mail: lzhu@hkbu.edu.hk